%%% Macros for Proceedings Volumes%%%%%%%%%%%%%%%%%%%%%%%%%%%%%%
%%% originally created June 1992 by Danny Lewis for %%%%%%%%%%%%
%%% Verlag Walter De Gruyter %%%%%%%%%%%%%%%%%%%%%%%%%%%%%%%%%%%%
%%% Revised version October 1995 %%%%%%%%%%%%%%%%%%%%%%%%%%%%%%
%%% Revised version June 2, 1996 %%%%%%%%%%%%%%%%%%%%%%%%%%%%%%
%%%PARAMETERS

\catcode`@=11
\nopagenumbers
\vbadness=10000
\hbadness=10000
\hsize=12.5cm
\vsize=19.2cm
\topskip=12pt
\parindent=0.5cm
\newskip\litemindent
\litemindent= 0,7 cm  %%%%% amount of indentation for \litem
\parskip=0pt
\widowpenalty=10000
\clubpenalty=10000
\hfuzz=1.5pt
\abovedisplayskip=6pt plus 1pt
\abovedisplayshortskip=0pt plus 1pt
\belowdisplayskip=6pt plus 1pt 
\belowdisplayshortskip=6pt plus 1pt
\frenchspacing

%%%%FONTS AND FAMILIES%%%%%

% AMS fonts are included, but commented out with
% the symbol "%.%", since not everyone has them.
% If you have them and need them,
% you can remove the %.%.

\font\authorfont=cmti10 at 14.40pt

\font\seventeenrm=cmbx10 at 17.28pt
\font\seventeenit=cmbxti10 at 17.28pt
\font\seventeeni=cmmib10 at 17.28pt
\font\seventeensy=cmbsy10 at 17.28pt
\font\seventeenex=cmex10 at 17.28pt
\font\seventeenmsa=msam10 at 17.28pt
\font\seventeenmsb=msbm10 at 17.28pt
\font\seventeeneuf=eufm10 at 17.28pt

\font\fourteenrm=cmbx10 at 14.40pt
\font\fourteenit=cmbxti10 at 14.40pt
\font\fourteeni=cmmib10 at 14.40pt
\font\fourteensy=cmbsy10 at 14.40pt
\font\fourteenex=cmex10 at 14.40pt
\font\fourteenmsa=msam10 at 14.40pt
\font\fourteenmsb=msbm10 at 14.40pt
\font\fourteeneuf=eufm10 at 14.40pt

\font\twelverm=cmbx10 at 12pt
\font\twelveit=cmbxti10 at 12pt
\font\twelvei=cmmib10 at 12pt
\font\twelvesy=cmbsy10 at 12pt
\font\twelveex=cmex10 at 12pt
\font\twelvemsa=msam10 at 12pt
\font\twelvemsb=msbm10 at 12pt
\font\twelveeuf=eufm10 at 12pt

\font\tenrm=cmr10
\font\tenit=cmti10
\font\tenbf=cmbx10
\font\tenib=cmmib10 
\font\teni=cmmi10
\font\tensy=cmsy10
\font\tenbsy=cmbsy10
\font\tenex=cmex10 
\font\tenmsa=msam10 
\font\tenmsb=msbm10 
\font\teneuf=eufm10

\font\ninerm=cmr9
\font\nineit=cmti9
\font\ninebf=cmbx9
\font\ninei=cmmi9
\font\ninesy=cmsy9

\font\ninemsa=msam9 
\font\ninemsb=msbm10 at 9pt
\font\nineeuf=eufm9

\font\eightrm=cmr8
\font\eightit=cmti8
\font\eightbf=cmbx8
\font\eighti=cmmi8
\font\eightib=cmmib8
\font\eightsy=cmsy8
\font\eightbsy=cmbsy8

\font\eightmsa=msam8
\font\eightmsb=msbm8
\font\eighteuf=eufm8

\font\sevenrm=cmr7

\font\sevenbf=cmbx7
\font\seveni=cmmi7

\font\sevensy=cmsy7

\font\sevenmsa=msam7
\font\sevenmsb=msbm7
\font\seveneuf=eufm7

\font\fiverm=cmr5
\font\fivebf=cmbx5
\font\fivesy=cmsy5

\font\fivei=cmmi5

\font\fivemsa=msam5
\font\fivemsb=msbm5
\font\fiveeuf=eufm5

\newfam\msafam%Family 8
\newfam\msbfam%Family 9 
\newfam\euffam%Family 10

\def\seventeenpoint{\def\rm{\fam0\seventeenrm}%
\textfont0=\seventeenrm\scriptfont0=\fourteenrm 
\scriptscriptfont0=\twelverm
\textfont1=\seventeeni\scriptfont1=\fourteeni 
\scriptscriptfont1=\twelvei
\def\mit{\fam1\seventeeni}\def\oldstyle{\fam1\seventeeni}%
\textfont2=\seventeensy\scriptfont2=\fourteensy 
\scriptscriptfont2=\twelvesy
\def\cal{\fam2\seventeensy}%
\textfont3=\seventeenex\scriptfont3=\fourteenex 
\scriptscriptfont3=\twelveex
\def\it{\fam\itfam\seventeenit}% 
\textfont\itfam=\seventeenit
\def\bf{\rm}
\def\msa{\fam\msafam\seventeenmsa}%
\textfont\msafam=\seventeenmsa\scriptfont\msafam=\fourteenmsa 
\scriptscriptfont\msafam=\twelvemsa
\def\msb{\fam\msbfam\seventeenmsb}%
\textfont\msbfam=\seventeenmsb\scriptfont\msbfam=\fourteenmsb 
\scriptscriptfont\msbfam=\twelvemsb
\def\euf{\fam\euffam\seventeeneuf}%
\textfont\euffam=\seventeeneuf\scriptfont\euffam=\fourteeneuf 
\scriptscriptfont\euffam=\twelveeuf
\setbox\strutbox=\hbox{\vrule height16pt depth4pt width0pt}%
\baselineskip=20pt\seventeenrm}

\def\fourteenpoint{\def\rm{\fam0\fourteenrm}%
\textfont0=\fourteenrm\scriptfont0=\twelverm 
\scriptscriptfont0=\tenbf
\textfont1=\fourteeni\scriptfont1=\twelvei 
\scriptscriptfont1=\tenib
\def\mit{\fam1\fourteeni}\def\oldstyle{\fam1\fourteeni}%
\textfont2=\fourteensy\scriptfont2=\twelvesy 
\scriptscriptfont2=\tenbsy
\def\cal{\fam2\fourteensy}%
\textfont3=\fourteenex\scriptfont3=\twelveex 
\scriptscriptfont3=\tenex
\def\it{\fam\itfam\fourteenit}% 
\textfont\itfam=\fourteenit
\def\bf{\rm}
\def\msa{\fam\msafam\fourteenmsa}%
\textfont\msafam=\fourteenmsa\scriptfont\msafam=\twelvemsa 
\scriptscriptfont\msafam=\tenmsa
\def\msb{\fam\msbfam\fourteenmsb}%
\textfont\msbfam=\fourteenmsb\scriptfont\msbfam=\twelvemsb 
\scriptscriptfont\msbfam=\tenmsb
\def\euf{\fam\euffam\fourteeneuf}%
\textfont\euffam=\fourteeneuf\scriptfont\euffam=\twelveeuf 
\scriptscriptfont\euffam=\teneuf
\setbox\strutbox=\hbox{\vrule height13pt depth4pt width0pt}%
\baselineskip=16pt\fourteenrm}

\def\twelvepoint{\def\rm{\fam0\twelverm}%
\textfont0=\twelverm\scriptfont0=\tenbf 
\scriptscriptfont0=\eightbf
\textfont1=\twelvei\scriptfont1=\tenib 
\scriptscriptfont1=\eightib
\def\mit{\fam1\twelvei}\def\oldstyle{\fam1\twelvei}%
\textfont2=\twelvesy\scriptfont2=\tenbsy 
\scriptscriptfont2=\eightbsy
\def\cal{\fam2\twelvesy}%
\textfont3=\twelveex\scriptfont3=\twelveex 
\scriptscriptfont3=\twelveex
\def\it{\fam\itfam\twelveit}% 
\textfont\itfam=\twelveit
\def\bf{\rm}%
\def\msa{\fam\msafam\twelvemsa}% 
\textfont\msafam=\twelvemsa \scriptfont\msafam=\tenmsa
\scriptscriptfont\msafam=\eightmsa%
\def\msb{\fam\msbfam\twelvemsb}% 
\textfont\msbfam=\twelvemsb \scriptfont\msbfam=\tenmsb%
\scriptscriptfont\msbfam=\eightmsb%
\def\euf{\fam\euffam\twelveeuf}% 
\textfont\euffam=\twelveeuf \scriptfont\euffam=\teneuf%
\scriptscriptfont\euffam=\eighteuf%
\setbox\strutbox=\hbox{\vrule height10pt depth4pt width0pt}%
\baselineskip=14pt\rm}

\def\tenpoint{\def\rm{\fam0\tenrm}%
\textfont0=\tenrm\scriptfont0=\sevenrm 
\scriptscriptfont0=\fiverm
\textfont1=\teni\scriptfont1=\seveni 
\scriptscriptfont1=\fivei
\def\mit{\fam1\teni}\def\oldstyle{\fam1\teni}%
\textfont2=\tensy\scriptfont2=\sevensy 
\scriptscriptfont2=\fivesy
\def\cal{\fam2\tensy}%
\textfont3=\tenex\scriptfont3=\tenex 
\scriptscriptfont3=\tenex
\def\it{\fam\itfam\tenit}% 
\textfont\itfam=\tenit
\def\bf{\fam\bffam\tenbf}% 
\textfont\bffam=\tenbf\scriptfont\bffam=\sevenbf
\scriptscriptfont\bffam=\fivebf
\def\msa{\fam\msafam\tenmsa}% 
\textfont\msafam=\tenmsa \scriptfont\msafam=\sevenmsa
\scriptscriptfont\msafam=\fivemsa
\def\msb{\fam\msbfam\tenmsb}%
\textfont\msbfam=\tenmsb \scriptfont\msbfam=\sevenmsb
\scriptscriptfont\msbfam=\fivemsb
\def\euf{\fam\euffam\teneuf}% 
\textfont\euffam=\teneuf \scriptfont\euffam=\seveneuf
\scriptscriptfont\euffam=\fiveeuf
\aline=12pt plus 1pt minus 1pt
\halfaline=6pt plus 1pt minus 1pt
\setbox\strutbox=\hbox{\vrule height8.5pt depth3.5pt width0pt}%
\baselineskip=12pt\rm}

\def\ninepoint{\def\rm{\fam0\ninerm}%
\textfont0=\ninerm \scriptfont0=\sevenrm \scriptscriptfont0=\fiverm
\textfont1=\ninei\scriptfont1=\seveni \scriptscriptfont1=\fivei
\def\mit{\fam1\ninei}\def\oldstyle{\fam1\ninei}%
\textfont2=\ninesy \scriptfont2=\sevensy \scriptscriptfont2=\fivesy
\def\cal{\fam2\ninesy}%
\textfont3=\tenex \scriptfont3=\tenex \scriptscriptfont3=\tenex
\def\it{\fam\itfam\nineit}%
\textfont\itfam=\nineit
\def\bf{\fam\bffam\ninebf}%
\textfont\bffam=\ninebf \scriptfont\bffam=\sevenbf
\scriptscriptfont\bffam=\fivebf
\def\msa{\fam\msafam\ninemsa}%
\textfont\msafam=\ninemsa \scriptfont\msafam=\sevenmsa
\scriptscriptfont\msafam=\fivemsa
\def\msb{\fam\msbfam\ninemsb}%
\textfont\msbfam=\ninemsb \scriptfont\msbfam=\sevenmsb
\scriptscriptfont\msbfam=\fivemsb
\def\euf{\fam\euffam\nineeuf}% 
\textfont\euffam=\nineeuf \scriptfont\euffam=\seveneuf
\scriptscriptfont\euffam=\fiveeuf%
\aline=11pt plus 1pt minus 1pt
\halfaline=5pt plus 1pt minus 1pt
\setbox\strutbox=\hbox{\vrule height7pt depth3pt width0pt}%
\baselineskip=11pt\rm}

\def\eightpoint{\def\rm{\fam0\eightrm}%
\textfont0=\eightrm \scriptfont0=\sevenrm \scriptscriptfont0=\fiverm
\textfont1=\eighti\scriptfont1=\seveni \scriptscriptfont1=\fivei
\def\mit{\fam1\eighti}\def\oldstyle{\fam1\eighti}%
\textfont2=\eightsy \scriptfont2=\sevensy \scriptscriptfont2=\fivesy
\def\cal{\fam2\eightsy}%
\textfont3=\tenex \scriptfont3=\tenex \scriptscriptfont3=\tenex
\def\it{\fam\itfam\eightit}% 
\textfont\itfam=\eightit
\def\bf{\fam\bffam\eightbf}%
\textfont\bffam=\eightbf \scriptfont\bffam=\sevenbf
\scriptscriptfont\bffam=\fivebf
\def\msa{\fam\msafam\eightmsa}%
\textfont\msafam=\eightmsa \scriptfont\msafam=\sevenmsa
\scriptscriptfont\msafam=\fivemsa
\def\msb{\fam\msbfam\eightmsb}%
\textfont\msbfam=\eightmsb \scriptfont\msbfam=\sevenmsb
\scriptscriptfont\msbfam=\fivemsb
\def\euf{\fam\euffam\eighteuf}% 
\textfont\euffam=\eighteuf \scriptfont\euffam=\seveneuf
\scriptscriptfont\euffam=\fiveeuf%
\aline=10pt plus 1pt minus 1pt
\halfaline=5pt plus 1pt minus 1pt
\setbox\strutbox=\hbox{\vrule height7pt depth3pt width0pt}%
\baselineskip=10pt\rm}

\skewchar\teni='177
\skewchar\ninei='177 
\skewchar\eighti='177 
\skewchar\seveni='177 
\skewchar\fivei='177
\skewchar\tensy='60
\skewchar\ninesy='60 
\skewchar\eightsy='60 
\skewchar\sevensy='60 
\skewchar\fivesy='60

%%%%headings

\long\def\Title#1\par{%
\global\titlepagetrue
{\parindent=0pt
\raggedcenter\pretolerance=10000
\seventeenpoint #1\par}
\vskip24pt}

\long\def\Author#1\par{%
\centerline{\fourteenpoint\authorfont #1}
\vskip60pt
\tenpoint\noindent\ignorereturn}

\def\Section{\removelastskip
\goodbreak\vskip36pt plus 1pt minus 6pt \section}

\def\section#1\par{%
{\raggedcenter
\interlinepenalty=10000\pretolerance=10000
\noindent\fourteenpoint #1\nobreak\par\nobreak}
\nobreak\vskip12pt\nobreak
\noindent\tenpoint\rm\ignorereturn}

\def\subsection#1\par{%
{\raggedcenter
\interlinepenalty=10000\pretolerance=10000
\noindent\twelvepoint #1\nobreak\par\nobreak}
\nobreak\vskip12pt\nobreak
\noindent\tenpoint\rm\ignorereturn}

\def\References #1 {\ifdim\lastskip<25pt \removelastskip
\vskip24pt plus 1pt\fi
\setbox0=\hbox{\ninepoint [#1]\enspace}%
\centerline{\twelverm References}\nobreak\par\nobreak
\interlinepenalty=10000\parskip=0pt plus 1pt
\litemindent=\wd0
\ninepoint
\nobreak\vskip7pt\nobreak}

%%%%GERMAN%%%%%

%%%%FOOTNOTES%%%%%

\newbox\footbox
\setbox\footbox=\hbox{\ninerm 22)~}

\newdimen\footnotespace
\footnotespace=0pt

\def\footnote#1{\let\@sf\empty %#2 (the text) is read later
  \ifhmode\edef\@sf{\spacefactor\the\spacefactor}\/\fi
$^{#1}$\@sf\vfootnote{#1}}

\def\vfootnote#1{\insert\footins\bgroup
\interlinepenalty\interfootnotelinepenalty
\splittopskip\ht\strutbox % top baseline for broken footnotes
\splitmaxdepth\dp\strutbox \floatingpenalty\@MM
\leftskip\z@skip \rightskip\z@skip \spaceskip\z@skip \xspaceskip\z@skip\parindent=\wd\footbox
\litem{\ninerm #1}\ninepoint\footstrut\futurelet\next\fo@t}
\def\fo@t{\ifcat\bgroup\noexpand\next \let\next\f@@t
\else\let\next\f@t\fi \next}
\def\f@@t{\bgroup\aftergroup\@foot\let\next}
\def\f@t#1{#1\@foot}
\def\@foot{\strut\egroup}
\def\footstrut{\vbox to\splittopskip{}}
\skip\footins=\bigskipamount % space added when footnote is present
\count\footins=1000 % footnote magnification factor (1 to 1)
\dimen\footins=8in % maximum footnotes per page

%%%%%OUTPUT ROUTINE%%%%%%

\newif\iftitlepage

\def\makeheadline{\iftitlepage \global\titlepagefalse
\vbox{\titleheadline}%
\vskip12pt
\else \vbox{\ifodd\pageno\rightheadline\else\leftheadline\fi}%
\vskip12pt\fi}

\def\lefthead{}
\def\righthead{}

\def\rightheadline{\line{\vbox to 8.5pt{}%
\hphantom{\tenrm\folio}\ninepoint\hfill\righthead\hfill{\tenrm\folio}}}
\def\leftheadline{\line{\vbox to 8.5pt{}%
\ninepoint{\tenrm\folio}\hfill\lefthead\hfill\hphantom{\tenrm\folio}}}
\def\titleheadline{\line{\vbox to 8.5pt{}%
\ninepoint\it\hfill[Page \folio]}}

\output{\plainoutput}

%%%%%%MACROS%%%%%%%

\def\Litem#1#2{\par\noindent\hangindent#1\litemindent
\hbox to #1\litemindent{\hfill\hbox to \litemindent
{#2 \hfill}}\ignorespaces}
\def\litem{\Litem1}

\newskip\aline \newskip\halfaline
\aline=12pt plus 1pt minus 1pt
\halfaline=6pt plus 1pt minus 1pt
\def\skipaline{\vskip\aline}

\def\ignorereturn{\def\neext{\afterassignment\restart
       \let\next}\neext}
\def\restart{\ifx\next\par\else\let\neext\next\fi\neext}

\def\raggedcenter{\leftskip=0pt plus 4em \rightskip
=\leftskip \parfillskip=0pt \spaceskip=.3333em
\xspaceskip=.5em \pretolerance=9999 \tolerance=9999
\hyphenpenalty=9999 \exhyphenpenalty=9999 }

\def\dotfill{\leaders\hbox to 1em{\hss.\hss}\hfill}

\def\Classification{\ifdim\lastskip<\aline\removelastskip\skipaline\fi
\noindent 1991 Mathematics Subject Classification: }

\long\def\Proclamation#1. #2\par{\ifdim\lastskip<\aline\removelastskip
\penalty-250 \skipaline\fi{\def\\##1){\litem{\rm(##1)}}\noindent
\bf#1\unskip. \it#2\par}\skipaline}

\def\Theorem {\Proclamation Theorem }
\def\Proposition {\Proclamation Proposition }
\def\Corollary {\Proclamation Corollary }
\def\Lemma {\Proclamation Lemma }

\def\Proof{\ifdim\lastskip<\aline\removelastskip\skipaline\fi
\noindent\it Proof. \rm}

\def\qedbox{$\rlap{$\sqcap$}\sqcup$}
\def\qed{\nobreak\hfill\penalty250 \hbox{}\nobreak\hfill\qedbox\skipaline}

\def\proclamation#1. {\ifdim\lastskip<\aline\removelastskip\penalty-250
      \skipaline\fi\noindent\bf#1\unskip. \rm}

\def\Remark{\proclamation Remark }

\def\ref #1 {\vskip4pt\litem{[#1]}}

\catcode`@=12
\def\@{\hbox{-}}

\def\Pn#1{{\bf P}^{#1}}
\def\Pnd#1{{\check{\bf P}}^{#1}}
\def\id#1#2{{\cal I}_{#1}(#2)}
\def\kn#1#2{{\cal O}_{#1}(#2)}

\def\coh#1#2#3{{\rm h}^{#1}(\kn {#2}{#3})}
\def\Coh#1#2#3{{\rm H}^{#1}(\kn {#2}{#3})}
\def\cohi#1#2#3{{\rm h}^{#1}(\id {#2}{#3})}

\def\rto{--->}

\def\C{{\bf C}}
\def\PP{{\bf P}}
\def\Z{{\bf Z}}

\Title Abelian surfaces with two plane cubic
 curve fibrations and Calabi-Yau threefolds

\Author Klaus Hulek and Kristian Ranestad*

\vfootnote {*}{Both authors were partially supported by
the HCM contract AGE (Algebraic Geometry in Europe), no ERBCHRXCT940557.
The first author would also like to thank MSRI for its hospitality.}

\vskip 10pt

1.  Introduction.

2.  Numerical possibilities.

3.  Elliptic scrolls and degenerate Calabi-Yau $3$-folds in $\Pn 5$.

4.  Non-normal Del Pezzo $3$-folds in $\Pn 5$.

5.  Non-normal Calabi-Yau $3$-folds in $\Pn 5$.

6.  Equations of elliptic scrolls in $\Pn 5$.

7.  Heisenberg symmetry of elliptic scrolls in $\Pn 5$.

8.  Conclusion.

\Section 1. Introduction

Abelian surfaces of small degree are contained in nodal Calabi-Yau $3$-folds,
simi\-larly many Calabi-Yau $3$-folds of small degree specialize to nodal
Calabi-Yau
$3$-folds with
abelian surfaces on them.  The first assertion is intimately connected with
the fact that
the moduli space of abelian surfaces of small degree is uniruled:  An
abelian surface on a
Calabi-Yau $3$-fold moves in a linear pencil, and therefore gives rise to a
$\Pn 1$ in the moduli
space of abelian surfaces.  This idea was taken up and explored by Gross
and Popescu [GP1, GP2]
starting with a very singular Calabi-Yau variety, the secant variety of an
elliptic normal
curve.  The translation scrolls inside the secant variety are degenerate
abelian surfaces
and form a $\Pn 1$ on the boundary of the moduli of abelian surfaces. They
show that the
secant variety deforms to nodal Calabi-Yau $3$-folds with only isolated
singularities and
with a pencil of abelian surfaces as long as the degree of the elliptic
curve is less than
11.  This limit is related to the Del Pezzo bound for the possible
smoothing of minimal
elliptic surface singularities.
\par
 We explore a similar setting.  In a
$\Pn 2$-scroll over an elliptic curve, any anticanonical divisor, if there
is one, is a,
possibly degenerate, abelian surface.  If one can glue two $\Pn 2$-scrolls
over elliptic
curves along an anticanonical divisor, the union is a singular Calabi-Yau
$3$-fold.  Furthermore, if the anticanonical divisor moves in a pencil on at
least one of the two
scrolls, then we are in a position like above.  \par
We start in section 2 by asking for smooth abelian surfaces with two pencils
of plane cubic
curves on it.  The two pencils would then define $\Pn 2$-scrolls whose union is
Calabi-Yau.  It turns out that purely numerical considerations bound the
degree of these
abelian surfaces by 18.  This bound is obtained by abelian surfaces which
form the complete
intersection $((0,3), (3,0))$ in the $\Pn 2\times \Pn 2$ with its Segre
embedding.
For each even degree $10\leq d\leq 18$ there are numerical
possibilities which are realized.   In this paper we study the associated
elliptic scrolls and Calabi-Yau $3$-folds in the case $d=12$, i.e. the case
of abelian surfaces embedded linearly normally in $\Pn 5$.
\par
In section 3 we find and describe the abelian surfaces of degree $12$ and the
two scrolls
defined by their pencils of plane cubic curves.  The union of the two
scrolls is a
non-normal Calabi-Yau $3$-fold of degree $12$.
\par
 In sections 4 and 5 a separate approach leads to constructions via
projected Del
Pezzo $3$-folds of non-normal Calabi-Yau $3$-folds in degrees 10,11, 12 and
13.
The projected  Del Pezzo $3$-folds are bilinked to the non-normal Calabi-Yau
$3$-folds. In the last three sections we prepare the argument that the
reducible Calabi-Yau $3$-folds described in section 3 may also be obtained via
bilinkage from projected Del Pezzo $3$-folds.
 \par
  After finding equations for elliptic scrolls in section 6, we devote
section
7 to the Heisenberg symmetry of elliptic scrolls and provide a
description of the family of Calabi-Yau $3$-folds that are unions
of two scrolls.
\par
More precisely, let $H_6$ be the Heisenberg group of level 6 and let $N_6$
be its
normalizer in $GL(6,{\bf C})$.  In $N_6$ there is a natural involution
$\iota$ which
restricts to the abelian surfaces as multiplication by $-1$.  Let
$G_6=\langle H_6,\iota\rangle$.
Then the space of cubics in $\Pn 5$ contains a 4-dimensional vector space of
$G_6$-invariant
pencils.\par
  We let ${\bf P}=\Pn 3$ be the parameter space for these
$G_6$-invariant pencils of cubics.
$H_6$ contains a subgroup isomorphic to $H_2$ and four subgroups
$H(K_1),\ldots, H(K_4)$ of index 3 containing this
subgroup.  For every subgroup $H(K_i)$ there is a
set of three lines in $\Pn 5$ containing $H(K_i)$ in its
stabilizer and left
invariant by the action of $H_6$.  Similarly there exists for each
subgroup $H(K_i)$ a
line $l_i$ in ${\bf P}$ parametrizing pencils of cubics which contain the
corresponding three lines.
\par
 A general point on any of these four lines $l_i$ in ${\bf P}$ corresponds to a
pencil of cubics which
defines an elliptic scroll singular along the corresponding three lines in
$\Pn 5$.  The scroll is residual to three $\Pn
3$'s in the complete intersection of the two cubics.
\par
Between each pair of lines $l_i,l_j$, there is a $1:1$
correspondence defined by the pairs of points which correspond to elliptic
scrolls that intersect along a possibly degenerate abelian surface.  The lines
spanned by corresponding points form a conic section in the Grassmannian of
lines in ${\bf P}$.  Altogether there are 6 disjoint conic
sections in the Grassmannian of lines in ${\bf P}$ which parametrize abelian
surfaces with two plane cubic curve fibrations. \par
In the final section we state and prove the main theorem of the paper:  Let
$X_E$ and $X_F$ be elliptic scrolls in $\Pn 5$ that intersect precisely along
a $(1,6)$-polarized abelian surface.

\Theorem 8.3.  The reducible $3$-fold $Y=X_E\cup X_F$
is a degeneration of irreducible non-normal Calabi-Yau $3$-folds of degree $12$
in $\Pn 5$. The general such $3$-fold is singular precisely along $6$ disjoint
lines.\par

Chang has described smooth $3$-folds of degree $12$ that are birational to
Calabi-Yau
$3$-folds (cf. [Ch], [DP]).  Similar to the ones described in this paper
they are
bilinked to Fano $3$-folds of degree $7$, but they differ by their
sectional genus.
This difference is manifested in the appearance of non-normal
singularities.  The
Calabi-Yau $3$-folds of Theorem 8.3 do not deform to smooth ones:  The
double point
formulas for $3$-folds in $\Pn 5$ give the class of the non-normal singular
locus in
terms of the coefficients of the Hilbert polynomial, so any deformation
also has
non-normal singularities.\par

 We find three open problems related to the
topics of this paper particularly interesting:
  \proclamation Problem 1.1. Find $H_6$-invariant non-normal but
irreducible Calabi-Yau $3$-folds, with a pencil of $(1,6)$-polarized abelian
surfaces, degenerating to the union of two elliptic scrolls.\par
 \proclamation Problem 1.2.  Consider the normalization of the Calabi-Yau
$3$-folds of degrees $10,\ldots,13$ of section 5. Find the invariants, the
Betti- and Hodge numbers and and describe the K\"ahler cone of these Calabi-Yau
$3$-folds.\par
 \proclamation Problem 1.3. Describe the elliptic scrolls and the reducible
Calabi-Yau $3$-folds in the cases $d=14,16,18$.\par
\skipaline
We work over the complex numbers.

\Section 2.  Numerical
possibilities

 Assume that an abelian surface $A\subset \Pn n$ has two fibrations
$$p:A\to E\quad  {\rm and}\quad  q:A\to
F$$
in plane cubic curves.  These fibrations define two $\Pn 2$-scrolls.  If
the fibers
of the maps $p$ and $q$ have intersection number $\geq 2$, then the planes
in the
two fibrations intersect in at least a line and the two $\Pn 2$-scrolls
coincide. So for our purposes we can assume that $F$ and $E$ are sections
of $p$
and $q$ with $E\cdot F=1$. In particular $A=E\times F$.
  \par
When there is no isogeny between $E$ and $F$ then
$A$ is the product of $E\times F$ in the Segre embedding of $\Pn 2\times
\Pn 2$ in $\Pn 8$. We shall now assume that $E$ is general in the sense
that $E$ has ${\rm End} E\cong \bf Z$ and that $\gamma:E\rightarrow F$ is a
primitive isogeny of degree $l$. Then $NS(A)$ is generated by $E, F$ and
$\Gamma$ where $\Gamma$ is the graph of $\gamma$. The numerical equivalence
 of a hyperplane divisor on $A$
may therefore be expressed as
$$H\equiv aE+bF+c\Gamma,\qquad a,b,c\subset\Z.$$
Notice that these surfaces and these divisors really exist for any general
elliptic curve $E$. We investigate for which numerical data they give us two
plane cubic curve fibrations. The intersection numbers are given by the table:
$$\matrix {&E&F&\Gamma\cr E&0&1&l\cr F&1&0&1\cr \Gamma &l&1&0.}$$ With the
requirements
$$H\cdot E=b+lc=3,\quad H\cdot F=a+c=3\quad {\rm and}
\quad H\cdot \Gamma=al+b\geq 2,$$
 we get $$d/2:=H^2/2=9-lc^2.$$ This means that $d/2\leq 9$ and that there
are the following possibilities for $H$:
 \skipaline
\centerline {Table 1.}
$$\vbox {\offinterlineskip\def\tablerule{\noalign{\hrule}}
\halign{\strut\quad $#$\quad &\quad\hfil $#$\hfil\quad &
\quad\hfil $#$\hfil \quad &\quad\hfil $#$\hfil \quad&\quad\hfil $#$\hfil  \quad
&\quad\hfil $#$\hfil  \cr\tablerule
 d&l&c&a&b&H\cdot\Gamma\cr
\tablerule
10&1&2&1&1&2\cr
10&4&1&2&-1&7\cr
12&3&1&2&0&6\cr
14&2&1&2&1&5\cr
16&1&1&2&2&4\cr
18&*&0&3&3&*\cr
\tablerule}}
$$
The $*$ in table 1 means that there is no isogeny between $E$  and $F$
involved. \skipaline

\Proposition 2.1. There exist abelian surfaces with two plane curve fibrations
of degree $d$ in $\Pn {{d\over 2}-1}$, when $d=10,12,14,16,18$.

\Proof To give examples
it is now enough to check that $H$ is very ample.  For this we use Reider's
criterion [Re], which in these cases reduces to check that there are no
elliptic curves $C$
on $A$ with $H\cdot C\leq 2$.  Any such curve $C$, not equivalent to $E$,
$F$ or $\Gamma$, must
intersect each of these strictly positively, i.e. $$C\cdot E>0,\quad
C\cdot F>0,\quad{\rm and}\quad C\cdot
\Gamma>0.$$ But in each case except the second $H=aE+bF+c\Gamma$ with
$a+b+c\geq 3$ so
$H\cdot C\geq 3$. In the second case $H=2E-F+\Gamma$, so if $C=\alpha
E+\beta F+\gamma
\Gamma$, then  $H\cdot C< 3$ implies that
$$3\alpha +3\beta +7\gamma\leq 2.$$
The other inequalities above yield
$$\beta +4\gamma \geq 1,\quad \alpha +\gamma \geq 1\quad{\rm and}\quad 4\alpha
+\beta \geq 1,$$
while $C^2=0$, since $C$ is elliptic, yields
$$\alpha\beta +4\alpha\gamma +\beta\gamma =0.$$
From the first four inequalities we get $\beta\leq -(\alpha +\gamma )$
which combines with the last equality to yield
$$0=4\alpha\gamma +\beta (\alpha +\gamma)\leq 4\alpha\gamma - (\alpha
+\gamma)^2=-(\alpha
-\gamma)^2.$$
This only occurs when $\alpha =\gamma$, i.e. from the relation $C^2=0$, when
$\alpha =0$ or $2\alpha +\beta=0$.  The former is impossible since $\alpha
+\gamma >0$, while
the latter is impossible since then $C\cdot H= 10\alpha
-6\alpha=4\alpha\geq 4$.
\par
 In the first case $H\cdot \Gamma=2$, so $H$ is not very ample. In fact
$|H|$ maps $A$ two to
one to a quintic elliptic scroll in $\Pn 4$.  In this case $E=F$ and the
scroll is  the
symmetric product of $E$.  Thus the two fibrations coincide in the image.
In each of the other
cases $|H|$ defines an embedding.\qed

\skipaline

\Remark 2.2. When $d=18$ there is the simple example of
$$E\times F=E\times \Pn
2\cap\Pn 2\times F\subset \Pn 2\times
\Pn 2$$
 in its Segre embedding in $\Pn 8$.  Clearly the union of the two
scrolls deform in this
case to Calabi-Yau $3$-folds.

   \skipaline
 In this paper we shall concentrate
on the case $d=12$.

\Section 3. Elliptic scrolls and
degenerate Calabi-Yau $3$-folds in $\Pn 5$

From now on we consider abelian surfaces $A\subset\Pn
5$ of degree $12$, i.e. with a $(1,6)$-polarization, and with two fibrations
 $$p:A\to E\quad {\rm and}\quad q:A\to F$$ in plane cubic curves.  Furthermore
we assume that $E$ is general.  We denote by $X_E$ the scroll of planes of the
fibration $p$, and by $X_F$ the scroll of planes of the fibration $q$.  The
corresponding $\Pn 2$-bundles are denoted $V_E$ and $V_F$, respectively.
\par
 In the notation of the
previous section and table 1  there is an isogeny $$\gamma:E\to F$$
 of degree $3$.  Furthermore, if $\Gamma$ is the graph of the isogeny, then
the hyperplane divisor is  $$H= 2E+\Gamma.$$
We fix an origin $o\in E$ and let $s_1=\gamma (t_1)$, for some $3$-torsion
point
$t_1$ on $E$, not in the kernel of $\gamma$.  Let $h:F\to F$ be translation by
$s_1$, and let $\gamma^*:{\rm Pic}^0F\to {\rm Pic}^0E$ be the isogeny dual to
$\gamma$.  If ${\cal O}_F(o)$ is the line bundle of degree $1$ whose unique
section vanishes at $o$, then $h^*({\cal O}_F(o))\otimes {\cal O}_F(-o)$
generates the kernel of $\gamma^*$, i.e. $${\rm ker}\,\,\gamma^*=\{{\cal O}_F,
h^*({\cal O}_F(o))\otimes {\cal O}_F(-o),{(h^2)}^*({\cal O}_F(o))\otimes {\cal
O}_F(-o)\}\quad (h^3={\rm id}).$$

\skipaline
\Proposition 3.1.  The rank $3$ vector bundle associated to the $\Pn
2$-bundle $V_F$ decomposes into
 ${\cal E}={\cal L}_0\oplus h^*{\cal L}_0\oplus (h^2)^*{\cal L}_0$,
where ${\cal
L}_0$ is a line bundle of degree $2$ on
$F$ and $h$ is as above.  Furthermore the scroll $X_F$ is singular precisely
along three lines, which span $\Pn 5$.

\Proof  First, notice that $\Gamma$ is contained in a pencil of hyperplanes
defined by the pencil $|H-\Gamma|=|2E|$ on $A$.  This pencil is the pullback by
$q$ of a divisor $\Delta_0$ of degree $2$ on $F$. Since any divisor in this
pencil
is a pair of plane cubic curves contained in a hyperplane, their planes
intersect in a point. Thus the linear system $|\Delta_0|$ defines a morphism of
$F$ of degree $2$ into the double locus of $X_F$.\par
Next we consider the
translates $\Gamma_t$ of $\Gamma$ on $A$ by a point on $t$ on $E$ and find
which
$\Gamma_t$ are contained in a pencil of hyperplanes. This happens precisely
when
$H-\Gamma_t$ is the pullback of a divisor of degree $2$ from $F$, or
equivalently
when the restriction of $\Gamma-\Gamma_t$ to a fiber $E_f=q^{-1}(f)$ of $q$ is
trivial. But the translate $\Gamma_t$ can be represented as $$\{(x+t,\gamma
(x))|x,t\in E\}$$ so the intersection
$$(\Gamma-\Gamma_t)\cap E_f=\gamma ^{-1}(f)-\gamma ^{-1}(f)+3t$$
which is trivial precisely when $3t=o$. Since $\gamma$ is already an
isogeny of degree $3$, the $3$-torsion points in the kernel of $\gamma$ leave
$\Gamma$ invariant under translation, so in fact we may choose $t=t_1$ as above
and we get 3 distinct translates of $\Gamma$ which are contained in a pencil of
hyperplanes in $\Pn 5$.  We denote them by
$\Gamma_0 (=\Gamma), \Gamma_1,\Gamma_2$.
They each span a $\Pn 3$ and determine a linear system of degree $2$ on
$F$, which we denote by $|\Delta_0|, |\Delta_1|, |\Delta_2|$. Each linear
system $|\Delta_i|$ defines a morphism of degree $2$ of $F$ to the double
locus of
$X_F$.  The images of these three maps are disjoint lines, and clearly the
planes of $X_F$ are spanned by the respective images by these three maps.
Therefore $V_F$ is defined by a decomposable rank $3$ bundle ${\cal E}$ of
degree
6.  If we denote the line bundle associated to the divisor $\Delta_i$ by ${\cal
L}_i$, then
 $${\cal E}={\cal L}_0\oplus {\cal L}_1\oplus {\cal L}_2.$$
To check the differences between the line bundles ${\cal L}_i$, we consider
the intersection
  $$(\Gamma-\Gamma_t)\cap F_e=\gamma (e)-\gamma (e-t)=\gamma
(t).$$ Since $3t_1=o$, translation by $s_1=\gamma (t_1)$ on $F$ is a
$3$-torsion
element that generates the kernel of the dual isogeny, $\gamma^*:{\rm
Pic}^0F\to {\rm Pic}^0E$, i.e.
   $${\cal L}_i=(h^{i})^*{\cal L}_0,$$
where $h$ is the translation on $F$ by $s_1$.

\par
It remains only to check the singularities of $X_F$.   Any
singular point of
$X_F$ is the intersection of two, possibly infinitely close, planes of
$X_F$.  But two planes
intersect only if they span at most a hyperplane.  The two planes are
defined by a section of
a line bundle $q^*{\cal L}$ for a line bundle ${\cal L}$ of degree $2$ on
$F$. They
intersect precisely when
$${\cal E}\otimes {\cal L}^{-1}$$
has a section. But this is the case precisely when ${\cal L}={\cal L}_i$ for
$i\in\{0,1,2\}$.  Furthermore,
for each of these three cases  ${\cal E}\otimes {\cal L}_i^{-1}$ has precisely
one section, so the
corresponding planes span a hyperplane and the intersection of the two
planes is only a
point.   \qed

The abelian surface $A$ is an anticanonical divisor on $V_F$.    We compute
the sections of $-K_{V_F}$.

\Lemma 3.2.  $\coh 0{V_F}{-K_{V_F}}=4$.\par
\Proof The natural isomorphism
$$\Coh 0{V_F}{-K_{V_F}}\cong {\rm H}^0(F,{\rm Sym}^3{\cal E}\otimes
{\cal L}_0^{-1}\otimes {\cal L}_1^{-1}\otimes {\cal L}_2^{-1}),$$
reduces the computation to counting trivial
summands of the rank $10$ vector bundle $${{\rm Sym}^3{\cal E}\otimes
{\cal L}_0^{-1}\otimes {\cal L}_1^{-1}\otimes {\cal L}_2^{-1}}.$$
 Since ${\cal L}_0^{-1}\otimes {\cal L}_1^{-1}\otimes {\cal L}_2^{-1}={\cal
L}_0^{-3}$ this count is the number of summands in

$${\rm Sym}^3{\cal E}={\rm Sym}^3({{\cal L}_0\oplus h^*{\cal L}_0\oplus
(h^2)^*{\cal L}_0})$$
which equal ${\cal L}_0^3$.  As $h^3={\rm id}$ this number is 4.\qed

\Lemma 3.3. The $\Pn 2$-bundle $V_F$ is the quotient of a trivial bundle $\Pn
2\times E$ by a cyclic group of order $3$.

\Proof Consider the isogeny $\gamma:E\to F$ again. Since $h^*({\cal
O}_F(o))\otimes {\cal O}_F(-o)$
generates the kernel of $\gamma^*$ the isogeny dual to $\gamma$,
$$\gamma^*{\cal L}_0\cong\gamma^*{h^*}{\cal L}_0\cong\gamma^*{(h^2)^*}{\cal
L}_0.$$
Therefore the pullback of $V_F={\bf P}({\cal E})$ over $F$ via
$\gamma$ trivializes the bundle.  The kernel of $\gamma$ is a cyclic group of
order three which acts on the pullback $\gamma^{\ast}{\cal E}$.\qed

Thus we may construct $V_F$ by starting with
$\Pn 2\times E$, and dividing by a suitable diagonal action of the
cyclic group of order $3$. For this we consider a vector space
$V=\langle x_0,x_1,x_2\rangle$ with the action
$$\tau: x_i\mapsto \epsilon^ix_i,\qquad i\in \Z_3.$$
As above, let $o\in E$ be the origin of $E$, and consider the
linear system $|3o|$ on $E$.   It embeds $E$ as a plane cubic
curve in
$\Pn 2$.  We may choose $\langle e_0,e_1,e_2\rangle$ as a basis for the
underlying
vector space of $\Pn 2$ such that
$$\tau: e_i\mapsto {\epsilon^{-i}}e_i$$
induces the action of translation by a $3$-torsion point $t_0$ on $E$.
The
 diagonal action defined by
$\tau\in \Z_3$ on $V\times E$:
$$\tau : v\times e\mapsto \tau(v)\times e+t_0$$
acts without fixed points, so the quotient
is a rank $3$ vector bundle on $E/\langle t_0\rangle=F$.
The  action of $\tau $ on $V$ decomposes into
the characters
$$V=\langle x_0\rangle\oplus\langle x_1\rangle\oplus\langle x_2\rangle.$$

The anticanonical divisors on $V_F$ pull back to anticanonical divisors
on $\Pn 2\times E$ which are invariant under the action of $\tau$.
But the anticanonical divisors on $\Pn 2\times E$ are just the
pullbacks of the cubic curves on the plane.
The action of $\tau$ on the plane has the following basis of invariant
cubics:
$$\langle x_0^3,x_1^3,x_2^3,x_0x_1x_2\rangle.$$
Since these have no basepoints, there are no basepoints for the system of
anticanonical divisors on $V_F$, and the general one is smooth.  Notice
furthermore that this
linear system of invariant cubics contains the Hesse pencil
$$\langle x_0^3+x_1^3+x_2^3,x_0x_1x_2\rangle,$$
and recall that the singular curves in this pencil are 4 triangles.  In
fact it is easy to
check that these four triangles are the only triangles in the linear system of
invariant
cubics.   The vertices of the triangle $x_0x_1x_2=0$ are mapped to the
singular lines of
$X_F$.  The vertices of the three other triangles sweep out elliptic normal
curves of degree
6 as we shall see next.

\Proposition 3.4.  The scroll $X_F$ is the $3$-torsion translation scroll of an
elliptic normal curve in $\Pn 5$.

\Proof  Consider an elliptic normal curve $C$ of degree $6$ in $\Pn
5$, embedded
by the linear system $|6o|$.   For any $P\in C$ consider the translation scroll
$$V_P=\cup_{y\in C}\langle y,y+P,y+2P\rangle.$$
This is, for general $P$, a $\Pn 2$-scroll of degree $18$.  When $3P=o$ the
points $y$ , $y+P$ and $y+2P$ generate the same plane, so then the translation
scroll has degree $6$.  In this case
$$\langle y,y+P,y+2P\rangle\cap\langle z,z+P,z+2P\rangle\not=\emptyset$$
precisely when $3y+3z=o$. But then the pencil $|y+z|$ defines a map
from $C$ to the double locus of $V_P$.
 Now $3(y+z)=3(y+z+P)=3(y+z+2P)$, so
this map factors through the isogeny $C\to C/\langle P\rangle$.  Thus the 9
linear
systems
$|y+z|$  with
$3(y+z)=0$ define three pencils of pairs of planes each defining a double
line for the translation scroll.  The translation scroll is clearly a
scroll over $C/\langle P\rangle$ and its associated vector bundle
decomposes into the
sum of three line bundles of degree $2$ which define the three double
lines. The
differences between any two of these line bundles are the powers of some
$3$-torsion line bundle. If $C\cong E$ and $C/\langle P\rangle=F$, we have
$V_P\cong
X_F$.\qed

\skipaline

Finally, we consider the union of the two elliptic scrolls $X_1=X_E$
and $X_2=X_F$ in $\Pn 5$ which both contain the $(1,6)$-polarized
abelian surface $A$. Let $V_1=V_E$ and $V_2=V_F$ be the corresponding
$\Pn 2$-bundles.  Then $V_1$ and $V_2$ are the normalizations of $X_1$
and $X_2$ and $A\in |-K_{V_1}|$, resp. $A\in |-K_{V_2}|$. Let
$Y=X_1\cup X_2$.  \skipaline

\Proposition 3.5. $Y$ has a partial desingularization $Y_0=V_1\cup
V_2$, which is a Calabi-Yau $3$-fold, i.e. $K_{Y_0}={\cal O}_{Y_0}$ and
$q={\rm h}^1({\cal O}_{Y_0})=0$. In particular $V_1\cap V_2=A$.

\Proof First, we may use the previous notation and let $A=E\times F$ with
hyperplane divisor  $H=2E+\Gamma$.\par
 Notice
that the abelian surface $A$ does not intersect the singular lines of
$X_F$.  In
fact $|H-E_f|$, for a plane cubic curve $E_f$ passing through a
singular point, would then have a basepoint. This is impossible by
Reider's criterion [Re].\par We show
next that $X_1$ and $X_2$ intersect transversally along $A$: Near $A$ the
$3$-fold
$Y$ is the intersection of two smooth $3$-folds.  If the intersection is not
transversal the tangent spaces of the two
scrolls at some point coincide.  This is possible only if the planes of the two
scrolls at the given point intersect along a line.  Those two planes
intersect $A$ in
elliptic cubic curves $E_f$ and $F_e$ which meet in a point.  Any curve in
$|H-E_f-F_e|$ has
degree $6$ and arithmetic genus $2$.  Since there is a pencil of
hyperplanes through
the two planes, this curve moves in a pencil on $A$,
this is impossible by Riemann Roch, so transversality follows.  If $X_1\cap
X_2$
contains a point disjoint form $A$, then a plane in $X_1$ meets a plane in
$X_2$
along a line and the argument above applies again.  Thus $X_1$ and $X_2$ meet
transversally along $A$. \par
Now, look at the partial desingularization $Y_0$ of $Y$ obtained by normalizing
$X_1$ and $X_2$. Thus we may write $Y_0=V_1\cup V_2$ and $V_1\cap V_2=A$.
To show that $q({Y_0})=0$ we consider ${\rm Pic}^0{Y_0}$.  An element in
${\rm
Pic}^0{Y_0}$ is a pair $({\cal M}_1, {\cal M}_2)$ with ${\cal M}_i\in {\rm
Pic}^0V_i$ which
glue along $A$.
Consider the fibrations
$$p_1:V_1\to F\quad{\rm and}\quad p_2:V_2\to E.$$
Now, ${\cal M}_i=p_i^*({\cal N}_i)$, where ${\cal N}_1\in {\rm Pic}^0(F)$
resp. ${\cal N}_2\in {\rm Pic}^0(E)$.  On $A$ we have sections of each
$p_i$ which are
fibres of the opposite map.  We identify these sections with $F$ and $E$
respectively.    Then ${\cal M}_1|_E={\cal O}_E$ and ${\cal M}_2|_F={\cal
O}_F$.  So
if we want to glue ${\cal M}_1$ and ${\cal M}_2$
 we must have that $${\cal O}_E={\cal M}_1|_E= {\cal N}_2$$ and
$${\cal O}_F={\cal M}_2|_F= {\cal N}_1,$$ and hence $({\cal M}_1, {\cal
M}_2)=({\cal
O}_{X_1}, {\cal O}_{X_2})$. This shows that   ${\rm Pic}^0{Y_0}=\{{\cal
O}_{Y_0}\}$.  Since ${\rm
Pic}^0{Y_0}$ is a reduced group scheme in characteristic 0 and since ${\rm
H}^1(X,{\cal O}_{Y_0})$
is the tangent space at the origin, it follows that  $q({Y_0})={\rm
h}^1({Y_0},{\cal O}_{Y_0})=0$. \par
  Since
$A\in |-K_{V_1}|$, resp. $A\in |-K_{V_2}|$ it is clear that the
restriction of $K_{Y_0}$ to $V_1$ and to $V_2$ is trivial.
But ${\rm Pic}^0{Y_0}=\{{\cal
O}_{Y_0}\}$, so $K_{Y_0}={\cal O}_{Y_0}$.\qed

\skipaline
The $3$-fold $Y=X_E\cup X_F$ is non-normal, singular along six lines three on
each scroll in addition to the surface $A=X_E\cap X_F$.  In the next section we
describe a series of non-normal Del Pezzo $3$-folds.  In section 5 we show
that these
are bilinked to non-normal Calabi-Yau $3$-folds. After some further
analysis of the
equations of $Y$ in sections 6 and 7 we show in section 8 that $Y$ is a
degeneration
of these non-normal Calabi-Yau $3$-folds.

\Section 4.  Non-normal Del Pezzo $3$-folds in $\Pn 5$

Del Pezzo $3$-folds are $3$-folds $W$ for which $K_W\cong -2H$, where $H$ is
ample.
Accordingly any smooth surface on $W$ in $|H|$ is a Del Pezzo
surface.\par
 Let
$V_t\subset \Pn {t+1}$ with $t=3,\dots ,8$ be the image of $\Pn 3$ by the map
defined by all
quadrics through $8-t$ general points in $\Pn 3$.  Then $V_t$ is a
Del-Pezzo $3$-fold.
 We describe the image $W_t\subset\Pn 5$ of the general projection of $V_t$.

In particular we want to describe the singular non-normal locus. Thus we are
interested in the cases $t=5,6,7,8$ and will prove

\Theorem 4.1. $W_t$ is non-normal along ${{t-3}\choose
2}$ skew lines and has ${{8-t}\choose
2}$ additional ordinary double points when $t=5,6,7,8$.\par

To prove this we will use a result of Reye on linear systems of quadrics.
To explain Reye's result we need the notion of apolarity applied to quadrics in
$\Pn 3$.  Thus let $S=k[x_0,\ldots ,x_3]$ and $T=k[y_0,\ldots, y_3]$, and
define a
pairing $S_2\times T_2\to k$ by letting $S$ operate as differential
operators on $T$
and vica versa.
\par
 We say that quadrics in $S$ and $T$ are apolar if they are
orthogonal
with respect to this pairing.  In our situation we think of $S$ as the
coordinate ring
of $\Pn 3$ and $T$ as the coordinate ring of $\Pnd 3$.  Starting with a
6-dimensional
subspace of quadrics $V\subset S_2$, there is a $4$-dimensional subspace
i.e. a web of
quadrics  $V^{\bot}\subset T_2$.
\par
 Since any quadric in $V$ is apolar to any
member of
$V^{\bot}$, we say that $V$ and $V^{\bot}$ are apolar sets of quadrics.
Consider the discriminant $D$ of the space $V^{\bot} $ of quadrics. This is
a quartic
surface defined by the determinant of a symmetric $4\times 4$ matrix with
linear entries.  The quadrics in
$V^{\bot}$ of rank $1$ and $2$ are respectively triple and double points on
$D$.  The
possible numbers of rank $1$ and rank $2$ quadrics are given in the
following table:
\par
\centerline {Table 2.}
$$\vbox {\offinterlineskip\def\tablerule{\noalign{\hrule}}
\halign{\strut\quad\hfil #\hfil\quad &\quad\hfil #\hfil\quad &
\quad\hfil #\hfil \quad \cr\tablerule
 {\rm rank $1$ quadrics}&{\rm
rank $2$ quadrics}\cr
\tablerule
  $4$&$0/\infty$\cr $3$&$1$\cr $2$&$3$\cr $1$&$6$\cr $0$&$10$\cr
\tablerule}}
$$
\skipaline
This follows from a few lemmas which have independent interest for us.

\Lemma 4.2. Each rank $1$ quadric in $V^{\bot}$ determines a basepoint
for the
quadrics in $V$ and vice versa.  Each rank $2$ quadric in $V^{\bot}$
determines a line
contained in 4 quadrics in $V$ and vica versa. Alternatively, if $\rho$ is
the map
defined by the linear system $V$ of quadrics, then each rank $2$ quadric in
$V^{\bot}$
determines a line in the source double point locus of $\rho$ and vice
versa.

\Proof Note that if $a=(a_0,a_1,a_2,a_3)$, $L=\sum_{i=0}^3a_iy_i$ and $q\in
S_2$, then $$L^2(q)=2q(a).$$
Now each rank $1$ quadric in $V^{\bot}$ has the form $L^2$ for some point
$a$, and
apolarity says that  $L^2(q)=2q(a)=0$ for every quadric $q$ in $V$ so $a$
is a basepoint
for $V$. Conversely, if $a=(a_0,a_1,a_2,a_3)$ is a basepoint for $V$ and
$L=\sum_{i=0}^3a_iy_i$, then $L^2(q)=2q(a)=0$ for every $q\in V$ so $L^2\in
V^{\bot}$.
\par
Each rank $2$ quadric in $V^{\bot}$ has the form $L_1^2+L_2^2$ for some
linear forms $L_1=\sum_{i=0}^3a_iy_i$ and
$L_2=\sum_{i=0}^3b_iy_i$.  Let $l\subset \Pn 3$ be the line spanned by the
points $a=(a_0,a_1,a_2,a_3)$
and $b=(b_0,b_1,b_2,b_3)$. Let $V_l\subset V$ be the subspace
of quadrics which vanish on $l$.  Then $$V_l=\{q\in V|
L_1^2(q)=2q(a)=L_2^2(q)=2q(b)=L_1L_2(q)=0\}.$$
But $(L_1^2+L_2^2)(q)=0$ for every $q\in V$, so $V_l$ has codimension 2 in
$V$.
Conversely if $V_l$ has codimension 2 in $V$, then some linear combination
of $L_1^2, L_2^2$ and $L_1L_2$
is contained in $V^{\bot}$.  But any such linear combination is a rank $2$
(or rank $1$)
quadric, and the lemma follows.\qed

 Porteous' formula [cf. Fu 14.4.11] computes the number of rank $2$
quadrics in a
general web of
rank $4$ quadrics.  This number is 10.  Reye found a geometric
interpretation  of these 10 rank $2$ quadrics considering apolar twisted
cubic curves to
the web of quadrics, i.e. twisted cubic curves whose defining net of quadrics
is apolar to $V^{\bot}$. \par
By a determinental net
of quadrics we mean a net (i.e. a 3-dimensional space) of quadrics which
is generated by the $2 \times 2$ minors of a $2 \times 3$ matrix of
linear forms.The general
determinental net of quadrics generates the ideal of a twisted
cubic curve.

\Lemma 4.3. {\rm (Reye).} The general 6-dimensional subspace
$V\subset S_2$ contains
precisely two determinental nets of quadrics, which together span $V$.  If
$C_1$ and
$C_2$ are the twisted cubic curves defined by these two nets, then every
rank $1$ quadric in $V^{\bot}$ corresponds to a point of intersection
between $C_1$ and
$C_2$ and vice versa.  Furthermore every rank $2$ quadric in $V^{\bot}$
corresponds to a
common secant line for $C_1$ and $C_2$ and vice versa.

\Proof  The number of determinental nets of quadrics in a general
6-dimensional space $V$ of quadrics is nowadays computable by quantum
cohomology [Kre] (compute the number of twisted cubic surface
scrolls through nine points in $\Pn 4$ in the quantum cohomology
of the Grassmannian of lines and intersect with a general
$\Pn 3$), a few years ago by modern intersection theory [ES] and in
the ancient times by direct geometric arguments [Rey].  We leave the
choice of reference to the reader.\par Given two determinental nets
which span $V$, the correspondence between points of intersection
and the rank $1$ quadrics follows from Lemma 4.2.\par
 For the second
correspondence consider first a common secant line to $C_1$ and $C_2$.
This line
and any one of the two curves form a complete intersection $(2,2)$.
Therefore the
line lies in a pencil of quadrics from each of the two determinental nets.
Together
the two pencils form a web of quadrics in $V$ which by Lemma 4.2
corresponds
to a rank $2$ quadric in $V^{\bot}$.  On the other hand the secant lines to a
twisted cubic curve form a congruence of bidegree $(1,3)$ in the Grassmannian
of lines in $\Pn 3$.  Thus two general twisted cubic curves have
$$(1,3)\cdot
(1,3)=1+9=10$$
common secant lines.  This is exactly the number of rank $2$ quadrics in the
web $V^{\bot}$, so the second correspondence follows.\qed

To fill in the remainder of table 2 we want to compute how much the
number
of rank $2$ quadrics decreases when the web acquires a rank $1$ quadric.  We
give an
argument using Reye's geometric interpretation.  Our arguments will
depend on a genericity assumption, i.e. the space $V$ is general with
a given number of basepoints.  The argument would go through without
this assumption also, but then the numbers in table 2 would have
different interpretations.  Since we will only use general systems
$V$ we do not consider the degenerate cases.\par
When the space $V$ of quadrics have
basepoints, then there will be infinitely many apolar twisted
cubic curves to $V$, but taking two of them will always suffice for our
argument.  First note that as long as the two twisted cubic curves have
less than 4 common points, the corresponding nets of quadrics do not
intersect. \par
 When the web contains
one rank $1$ quadric, then the two twisted cubic curves have one common
point.  The number
of common secant lines passing through this point is, with our
genericity assumption, easily computed by projection from
the point, it is 4, the number of intersection points between two conics in
the plane.
So the web has 6 rank $2$ quadrics in addition to the rank $1$ quadric.
This is the second row in the table.\par When the web contains 2 rank $1$
quadrics,
then the two twisted cubic curves have two common points.  In this case
there are 4
common secant lines through each of the two intersection points, and one
of these
is the line passing through the two points, so there are exactly 3 common
secant lines
which do not pass through any of the two common points.  Thus the web has 3
rank $2$
quadrics in addition to the two rank $1$ quadrics. This is the third row in
the table.
\par
When the web contains 3 rank $1$ quadrics, the two twisted cubic curves have
3 points in
common.  There are 4 common secant lines through each intersection point,
and adding
up three are counted twice, so we get only one common secant line which
does not pass
through any of the intersection points. Thus the web has one rank $2$
quadric in addition to the three rank $1$ quadrics. This is the fourth row in
the table.
\par
When the web has more than 3 rank $1$ quadrics, the twisted cubics have
at least 4
points of intersection.  In this case the number of common secant lines
that does not
pass through the intersection points is infinite or zero depending on
whether the
two determinental nets intersect or not. This covers the remaining row
in the table.

\Lemma 4.4. Consider the map $\rho$
defined by the linear system $V$ of quadrics, and consider a subscheme $Z$
of length 2
which does not intersect the baselocus of $V$. Then $Z$ is mapped to a
point by $\rho$
if and only if either the restriction of $\rho$ to the  unique line passing
through $Z$ is
$2:1$ onto a line, or this line contains two base points for $V$.

\Proof  The linear
system $V$ restricted to the line through $Z$ has degree $2$. If this linear
system has
no basepoints, then  the image of the line by $\rho$ is a line or a conic
section. It
is a line if and only if some subscheme of length 2 is mapped to a point.
If the line
through $Z$ intersects the baselocus, the intersection must contain two
points, such
that the line is contracted by $V$.\qed

{\it Proof of Theorem 4.1.}
Combining Lemmas 4.2, 4.3 and 4.4 we find that a general 6-dimensional linear
system of quadrics with $8-t$
basepoints and $t=5,6,7,8$ defines a rational map from $\Pn 3$ to $\Pn
5$, whose image has
degree $t$ and whose non-normal double locus consists of ${{t-3}\choose
2}$ disjoint lines.  Any line between basepoints is contracted to an isolated
singularity.   The number of isolated singular points is therefore
${{8-t}\choose
2}$.  It is easily checked by restriction to general hyperplane sections
that these singular points are ordinary double points.\qed

\skipaline
By abuse of notation we call the varieties $W_t$ Del Pezzo $3$-folds.  Their
normalization have only isolated double points from the contracted lines
between
basepoints, these occur when $t\leq 6$. For each $t$ we want to describe
the ideal of
$W_t$ and understand their linkage class.

\Proposition 4.5.  A $3$-fold $W_5$ is bilinked to the union of two
$\Pn 3$'s which span $\Pn 5$. A $3$-fold $W_6$ is bilinked to the union of
three $\Pn 3$'s.
The ideals of general Del Pezzo $3$-folds $W_7$ and $W_8$ are generated by
quartics and
quintics, and quintics and sextics respectively.\par

\Proof To understand the ideal of $W_t$ we first investigate the ideal
of the singular
lines.  For $t=5$ there is one line so this case is trivial.
 The case $t=6$ is
also easy since there are three lines; they span $\Pn 5$ as soon as the
projection is general, as is easily verified.  For $t=7$ and $t=8$ the
situation
is a bit more involved.

\Lemma 4.6.  The 6 singular lines of a general Del Pezzo $3$-fold $W_7$
lie in a determinental net of quadrics in $\Pn 5$, in fact in the Segre
embedding of $\Pn 1\times \Pn 2$.

\Proof  Let $V$ be the linear system of quadrics which define the
rational map of
$\Pn 3$ onto $W_7$.  By Lemma 4.3 the linear system $V$ contains and is spanned
by at least two determinental
nets. Of course, the corresponding twisted cubic curves pass through
the basepoint, so the union of any two is
 rational and is therefore contained in a cubic surface.  By
genericity we may assume that this surface is smooth.  On this cubic
surface there is a pencil of twisted cubic curves through the basepoint
with six common secant lines. This pencil of curves corresponds to a linear
pencil of
nets of quadrics, so since two nets are contained in $V$ they all are. The
images of
these twisted cubic curves on $W_7$
 are plane curves.  Since the twisted cubic curves pass through the
basepoint of
$V$ these plane curves have degree $5$.  The 6 common secant lines are
mapped to
the 6 singular lines on $W_7$ and they account for the 6 singular points on
each of these
plane rational quintic curves.
  \par Consider the space of quadrics passing through the 6
singular lines.  If any of these quadrics intersects a plane of the plane
quintic curves
properly, its intersection would be a conic section passing through the 6
singular points of
the quintic curve. But this is impossible by Bezout.  Therefore any
quadric which passes
through the 6 lines must contain the planes of these curves. The
intersection of these
quadrics is therefore at least a threefold.  Now, 6 lines impose at most 18
conditions on
quadrics, so there are at least 3 such quadrics and the intersection is a
threefold.  A
codimension 2 variety contained in 3 quadrics is a rational cubic scroll.
If this scroll
is singular, it is a cone, and any two planes meet.  In this case two
planes span a
hyperplane.  Pulled back to the cubic surface of the twisted cubic curves,
this hyperplane
corresponds to a quadric which contains the two curves.  But the union of
the two
curves has arithmetic genus $0$, so it is not contained in a quadric and a
cubic
surface. Therefore the scroll is smooth, isomorphic to the Segre $3$-fold
scroll.
\qed

\Proposition 4.7.  $W_7$ is contained in five quartic hypersurfaces, they
define
an arithmetically
Cohen Macaulay $3$-fold, the union of $W_7$ and the rational cubic scroll $R$
which contains
the singular lines of $W_7$.  The quartic hypersurfaces are defined by the
maximal minors of a $4\times 5$ matrix with
linear entries.

\Proof  Again let $V$ be the linear system of quadrics
which define the
rational map of $\Pn 3$ onto $W_7$. In the proof of Lemma 4.6 we saw that
there is a pencil
of twisted cubic curves through the basepoint $p$ of $V$ whose defining nets of
quadrics are contained
in $V$.  These twisted cubic curves sweep out a cubic surface $S_3$ whose image
contains the singular lines on $W_7$ and is contained in the rational cubic
scroll $R$. The linear system of quadrics $V$ restricts to $S_3$ with one base
point, so the image has degree $11$. Consider the quartic surfaces through
$S_3$ and
singular in $p$; they consist of the union of $S_3$ and planes through
$p$, so these quartics form a net.  On the other hand, the image of $S_3$ in
$\Pn 5$ is contained in the net of quadrics through $R$,  and these
quadrics pulled back to
$\Pn 3$, correspond to quartic surfaces through $S_3$ singular at $p$, i.e.
to the
above net of quartics. Since the net of quartics has no unassigned
basepoints, the
quadrics through $R$ define precisely the image of $S_3$ on $W_7$, i.e. the
intersection $W_7\cap R$ is precisely the image of the surface $S_3$.  The
union
$W_7\cup R$ has degree $10$ and genus 11.  To conclude
that the union is arithmetically Cohen Macaulay we give an example.
Let $V$ be the space of quadrics
$$
\langle x_0^2+x_1x_2,
x_1^2+x_2x_3,
x_2^2+x_3x_0,
x_0x_1,
x_0x_1+x_2x_3,
x_0x_2+x_1x_3\rangle.$$
A straightforward computation in [MAC] shows that the Del Pezzo $3$-fold
$W_7$ defined by
$V$ lies in precisely 5 quartics, the $4\times 4$ minors of a $4\times 5$
matrix with linear
entries, i.e. these quartics define an arithmetically Cohen
Macaulay scheme of degree $10$ and genus $11$. Since this is an open condition
in the Hilbert
scheme  [Ell], the same is true for the general projection
$W_7$.\qed

For $W_8$ we get somewhat less.

\Lemma 4.8.    $W_8$ contains two plane rational
sextic curves, and each of the singular lines is spanned by a pair of nodes
of these
two sextic curves.  In particular there are
sextic generators in the ideal of $W_8$.

\Proof  In this case there are two apolar twisted cubic curves, these are
mapped to plane
sextic curves with 10 double points at the intersection of these planes
with the 10
singular lines of $W_8$.    So
we need sextic generators in the ideal of $W_8$.\qed

We are now ready to give some numerical results for the ideals of general
Del Pezzo $3$-folds
$W_t$ of degree $5\leq t\leq 8$.
\Lemma 4.9.  Table $3$ gives the degrees $d$ and the
number of generators in the ideal
of $W_t,$ for $t=5,6,7,8$.\par
\centerline {Table 3.}
$$\vbox {\offinterlineskip\def\tablerule{\noalign{\hrule}}
\halign{\strut\quad $#$\quad &\quad\hfil $#$\hfil\quad &
\quad\hfil $#$\hfil \quad &\quad\hfil $#$\hfil \quad&\quad\hfil
$#$\hfil  \quad
&\quad\hfil $#$\hfil  \quad \cr\tablerule
&d=2&d=3&d=4&d=5&d=6\cr
\tablerule
t=5&0&5&0&0&0\cr
t=6&0&1&7&0&0\cr
t=7&0&0&5&5&0\cr
t=8&0&0&1&10&\geq 1\cr
\tablerule}}$$
\Proof  The following spaces $V$ of quadrics
$$
\langle x_1x_2+x_0x_3,
x_3^2+x_2x_0,
x_3^2+x_1x_0,
x_3^2+x_2x_1,
x_0x_1+x_2x_3,
x_0x_2+x_1x_3\rangle,$$
$$
\langle x_0^2+x_1x_2,
x_1^2+x_2x_3,
x_3x_0,
x_0x_1,
x_0x_1+x_2x_3,
x_0x_2+x_1x_3\rangle,$$
$$
\langle x_0^2+x_1x_2,
x_1^2+x_2x_3,
x_2^2+x_3x_0,
x_0x_1,
x_0x_1+x_2x_3,
x_0x_2+x_1x_3\rangle,$$
and
$$
\langle x_0^2+x_1x_2,
x_1^2+x_2x_3,
x_2^2+x_3x_0,
x_3^2+x_0x_1,
x_0x_1+x_2x_3,
x_0x_2+x_1x_3\rangle,$$
have respectively 3,2,1 and no basepoints. The ideals of the corresponding
Del Pezzo $3$-folds
in $\Pn 5$ are easily computed in [MAC] and have Betti
numbers as in the table.
\par
For the proof of the lemma we show that the table
represents a lower bound on the number of generators in each of the given
degrees. The
lemma then follows by semicontinuity.
\par
In each case we define $\Sigma$ to be the union of the singular lines $L_i$ on
$W_t$, i.e. precisely the non-normal double point locus of $W_t$.  Let $V_t$ be
the normalization of $W_t$, i.e. $V_t$ is isomorphic to $\Pn 3$ blown up in
$8-t$
points, and with ${{8-t}\choose 2}$ lines contrated to ordinary double points.
The map $\varphi: V_t\to W_t$ is double precisely along $\Sigma$, thus we
get an exact sequence of sheaves
$$0\to {\cal O}_{W_t}\to \varphi_{\ast}{\cal O}_{V_t}\to
{\cal O}_{\Sigma}(-1)\to 0,$$
where ${\cal O}_{\Sigma}=\oplus_{i=1}^l{\cal O}_{L_i}$ and $l={{t-3}\choose
2}$.  Since
$V_t$ is a Del Pezzo $3$-fold, the cohomology of this sequence gives
$\coh 3{W_t}k={\rm h}^3(\varphi_{\ast}{\cal O}_{V_t}(k))=0,$ when $k\geq
3$.  Furthermore
the Euler characteristic of the relevant twists of the ideal sheaf of $W_t$
is easily
computed from this exact sequence together with the exact sequence $$0\to
\id {W_t}k\to \kn
{\Pn 5}k\to \kn {W_t}k\to 0.$$
  We collect the
results in the following table.
\skipaline
\centerline { Table 4.}
$$\vbox {\offinterlineskip\def\tablerule{\noalign{\hrule}}
\halign{\strut\hskip 4pt $#$\hskip 4pt &\hskip 4pt\hfil $#$\hfil\hskip 4pt &
\hskip 4pt\hfil $#$\hfil \hskip 4pt &\hskip 4pt\hfil $#$\hfil\hskip
4pt&\hskip 4pt\hfil $#$\hfil \hskip 4pt &\hskip 4pt\hfil $#$\hfil \hskip 4pt
&\hskip 4pt\hfil
$#$\hfil\cr\tablerule &k&\chi ({\cal O}_{W_t}(k))&{\rm h}^0(\varphi_{\ast}{\cal
O}_{V_t}(k))&h^0 ({\cal
O}_{\Sigma}(k-1)) &\coh 0{\Pn 5}k&
\chi ({\id {W_t}k})\cr
\tablerule
t=5&&&&&&\cr
&3&51&54&3&56&5\cr
t=6&&&&&&\cr
&3&55&64&9&56&1\cr
&4&113&125&12&126&13\cr
t=7&&&&&&\cr
&4&121&145&24&126&5\cr
&5&221&251&30&252&31\cr
t=8&&&&&&\cr
&4&125&165&40&126&1\cr
&5&236&286&50&252&16\cr
\tablerule}}$$

By restriction to general $\Pn 3$ sections of $W_t$ it is easy to check
that $\coh
1{W_t}k=0$ for the values of $k$ in the table.  Since additionally
$\coh 3{W_t}k=0$ when $k\geq 3$, the Euler characteristic of
the twisted ideal is a lower bound  for $\cohi 0{W_t}k$.  Thus the Betti
numbers of the
ideal
follow except in the case of quintics and $W_7$ and sextics and $W_8$.
These cases are
accounted for in Proposition 4.7 and Lemma 4.8.   In fact, from Lemma 4.7
we get 4
linear syzygies among the 5 quartics in
the ideal of $W_7$. Therefore there are also 4 extra quintic generators in the
ideal, 5 altogether.
Lemma 4.8 says that there are sextic generators in the ideal of $W_8$.
 \qed

To finish the proof of Proposition 4.5 it remains to consider the linkage
classes of $W_5$ and $W_6$.  The $ 3$-fold  $W_5$ is linked $(3,3)$ to a
rational
$3$-fold scroll of degree $4$. This lies in a quadric and is linked $(2,3)$ to
two $\Pn 3$'s, which clearly is minimal in its even biliaison class.\par
$W_6$ is linked $(3,4)$ to a $3$-fold $U$ with sectional genus $1$ which lies
on two cubics.
Consider a $\Pn 3$ spanned by two singular lines in $W_6$. It intersects
$W_6$ in the two
lines and in two additional skew lines each intersecting both the singular
lines.  Clearly
every cubic through $W_6$ contain this $\Pn 3$, so $U$ must intersect it in
a quartic surface
singular along the two singular lines.  $U$ is therefore an elliptic scroll
linked $(3,3)$ to
the union of three $\Pn 3$'s  (cf. also section 6).  This concludes
the proof of Proposition 4.5. \qed

In the last section we need a converse to Proposition 4.7.
\Proposition 4.10.  Let $W$ be a $3$-fold with sectional genus $1$, with
non-normal
double points along $6$ skew lines, no $3$ in a $\Pn 3$ and not all $6$ on a
rational normal quartic scroll. Assume that there is a Segre cubic scroll
$R$ in
$\Pn 5$ containing the $6$ skew lines transverse to its planes and
intersecting $W$
in a surface of degree $11$, linked to $4$ planes in the intersection of
$R$ with a
quintic hypersurface.  Assume furthermore that $W\cup R$ is
scheme-theoretically
defined by the $4\times 4$ minors of a $4\times 5$ matrix with linear
entries, and
that the only common singularities of these quartics are the $6$ lines.
Then $W$ is
a Del Pezzo $3$-fold $W_7$. \par
 \Proof  Consider the $4\times 5$ matrix $M$ whose maximal minors define
$R\cup W$.  Let $\langle z_0,\ldots, z_3\rangle$ be the coordinates of $\Pn
3$, then
the 5 bilinear equations $$(z_0,z_1,z_2,z_3)\cdot M=(0,0,0,0,0)$$ define a
$3$-fold
$T$ in $\Pn 3\times \Pn 5$.  The projection to $\Pn 5$
is onto $W\cup R$, while the other projection is onto $\Pn
3$.  The fiber in $T$ over any point in $\Pn 3$ is linear, defined by
the $5$ linear equations in the $6$ coordinate functions of $\Pn 5$, so for
the general point in $\Pn 3$
the fiber is a single point. The fiber over any point in $W\cup R$ is also
linear, with
dimension equal to $3$ minus the rank of $M$ at the target point.  The
points where $M$ has rank at most $2$, are singular on $W\cup R$, in fact
they are
singular on any quartic minor of $M$, so by our assumption, these points
all lie on
the $6$ singular lines.  Outside the $6$ singular lines $M$ has rank $3$
and therefore there is  a rational map
$$\varphi: W\cup R\rto \Pn 3$$
which is a morphism outside these $6$ lines.  This map is birational on one
component
of $W\cup R$ and contracts every other component to a surface, a curve or a
point.
The inverse map
  $$\psi:\Pn 3\rto W\cup R$$
is defined by the maximal minors of a $5\times 6$-matrix $M'$ with linear
entries.  $M'$
is obtained from $M$ by interchanging rows and linear forms.   Surfaces
that are images of components of
$W\cup R$ are fixed components for the linear system of quintic minors of
$M'$.\par
Our first aim is to show that $W$ has to be irreducible.
We analyse carefully possible reducible components of the surface $S=W\cap
R$.  For this we start with
the intersection of $S$ with a general plane in $R$.\par
  On $R$ the surface $S\equiv 5h-4f$, where $h$ is
the class of a hyperplane section, while $f$ is the class of a plane.
Since $W$ is non-normal along $6$
lines transverse to the planes of $R$, the intersection of $S$ with a
general plane is a curve $C$ of
degree $5$ singular in $6$ points, the points of intersection between the
plane and the singular lines.
Furthermore, no three singular lines span a $\Pn 3$, so no three of the
singular
points of $C$ are collinear.  More generally, any effective divisor of type
$ah-af$ is the product
of a plane curve of degree $a$ and $\Pn 1$. Geometrically each $\Pn 1$ is a
line
transverse to all planes in $R$. Conversely, any irreducible divisor of
type $h+bf$ that contains
more than one line transverse to all planes has $b=-1$ and is a quadric
surface, while any
irreducible divisor of type $2h+bf$ that contains more than $4$ lines
transverse to all planes
has $b=-2$ and is a rational normal quartic scroll. In our situation, this
means that three
singular lines intersect no plane in collinear points.  Also, since not all $6$
singular lines are on a rational normal quartic scroll, the $6$ singular
points cannot lie on a
conic section. Therefore, if $C$ is reducible it is either the union of $2$
conics and a line, or
the union of a conic and a cubic.  In the former case the singular points
are the points of
intersection between the two conics and one point of intersection between
the line and each conic,
while in the latter case $5$ singular points are points of intersection
between the conic and the
cubic, while the last point is singular on the cubic.  \par Next, the
surface $S$ intersects each
general line in $R$ transverse to the planes in one point. Since the
intersection of $S$ with a
general plane has at most $3$ components, $S$ itself has at most 4
components: Only the planes in
$R$ do not intersect the general plane in a curve, so since exactly one
component intersects every
line transverse to the planes, there are at most $4$ components of $S$.\par
 To describe the possible components of $S$ we
first note that on $R$ any effective divisor is of type $ah-bf$, with
$a\geq 0$ and $b\leq a$.
 In our two cases of possible reducible curves in general planes
we have the following possible decompositions of $S$ into irreducible
components:
$$S\equiv
5h-4f=b_0f+(h-b_1f)+(2h-b_2f)+(2h-b_2^{\prime}f)$$ where
$b_1+b_2+b_2^{\prime}-b_0=4$ and
$$S\equiv 5h-4f=b_0f+(2h-b_2f)+(3h-b_3f)$$ where $b_2+b_3-b_0=4$.  With the
above restriction on
the $b_i$, the former case occurs only when
$(b_0,b_1,b_2,b_2^{\prime})=(1,1,2,2)$,
$(b_0,b_1,b_2,b_2^{\prime})=(0,0,2,2)$ or
$(b_0,b_1,b_2,b_2^{\prime})=(0,1,2,1)$.  The latter
case occurs only when $(b_0,b_2,b_3)=(1,2,3)$, $(b_0,b_2,b_3)=(0,2,2)$ or
$(b_0,b_2,b_3)=(0,1,3)$.\par

We are now ready to analyse possible components of $W$.
 Clearly any two components of $W$ intersect each
other along a surface or a curve (or both).  If two
components of $W\cup R$ intersect each other along a curve, then every
hypersurface through $W\cup R$ has to be singular along this curve, so by our
assumption this happens only along the $6$ singular lines. The  sectional
genus of
$W$ is $1$, thus for the general $\Pn 3$ section of $W$ the arithmetic
genus is 1
plus the number of intersection points with the singular lines.  This means
that
for a general $\Pn 3$ passing through a plane of $R$, the contribution to the
arithmetic genus of the plane curve $C$ is $0$.  The degree of $W$ is $7$, so
residual to $C$ in this intersection is a curve of degree $2$.  It must
intersect
$C$ in 2 points, and it must itself have arithmetic genus 0, so it is a conic
section, two connected lines or a double line.  This very much restricts the
possiblities for the irreducible components of $W$.\par  In our analysis of the
components $W_d$ of $W$, we shall index them by their degree, unless the
degree is
not specified, in which case we use index $0$. This notation applies only
in this
proof and should not be confused with Del-Pezzo $3$-folds $W_t$. A
component $W_0$
of $W$ must intersect $R$ in a curve or a surface.  If it intersects $R$
only in a
curve, then this curve meets the general plane in $R$ only in points, thus this
component intersects a general $\Pn 3$ through a general plane in $R$ in a
curve
not contained in the plane.  Therefore this component has degree $1$ or $2$.
Similarly if $W_0$ intersects the general plane in $R$ in a curve of degree
$a$,
then the degree of the component is $a$, $a+1$ or $a+2$.  In fact, if $W_0$
spans
$\Pn 5$, i.e. when $a\geq 2$, then the degree is  strictly greater than
$a$.   Of
course, while $W_0$ is irreducible, $W_0\cap R$ may still be reducible.
Anyway,
we may enumerate the different cases in a table, where the columns are
ordered by
the degree of the curve of intersection between a component and a general
plane,
while the entries are the degrees of the components themselves:

\centerline {Table 5.}
$$\vbox {\offinterlineskip\def\tablerule{\noalign{\hrule}}
\halign{\strut\quad $#$\quad &\quad\hfil $#$\hfil\quad &
\quad\hfil $#$\hfil \quad &\quad\hfil $#$\hfil \quad&\quad\hfil $#$\hfil  \quad
&\quad\hfil $#$\hfil  \cr\tablerule
 0&1&2&3&4&5\cr
\tablerule
1&0&0&0&0&6\cr
1&1&0&0&5&0\cr
0&1&0&0&0&6\cr
0&1&3,3&0&0&0\cr
0&0&3&4&0&0\cr
0&0&0&0&0&7\cr
\tablerule}}
$$

 In the first two cases of table 5 the first linear component $W_1$ may
intersect $R$
in a plane and a line, or in a connected curve of degree $3$.  The latter
may be
excluded since no connected cubic curve is supported on the $6$ singular
lines.  In
the former case the line must be one of the $6$ singular lines.  The
remaining part
of $W$ is a component of degree $6$ or two components of degree $1$ and $5$,
respectively.   These remaining components  must intersect $R$ in a surface
of type
$5h-5f$, i.e. they do not intersect the general line in $R$ transverse to the
planes.  In the first case the linear component intersects the component
$W_6$ of
degree $6$ along a quadric surface, since the sectional genus is $1$.  The
linear
component $W_1$ intersects $R$ in a plane, but $W_6$ intersects the general
plane in
$R$ in a curve of degree $5$, so this plane must also be contained in
$W_6$.  This
contradicts the fact that $W_6$ does not intersect the general line in $R$
transverse
to the planes. The second case is impossible for the same reason.\par In
the third
case of the table there is a linear component $W_1$ that intersects $R$ along a
quadric surface and a component $W_6$ that intersects $R$ in a surface of type
$4h-3f$.  Again the two components have also to intersect along a quadric
surface.
The intersection of the two components on $R$ is a curve of type
$(h-f)(4h-3f)=4h^2-7hf$, i.e. a curve of degree $5$.  On the other hand on
$W_1$ each
of the two other components intersect along a  quadric surface and two singular
lines.  This adds up to a curve of degree $6$ in the intersection of all three
components.  This is a contradiction. The fourth case is entirely
similar.\par In the
fifth case there are two components $W_3$ and $W_4$ of degree $3$ and $4$
respectively.  They intersect each other in a surface of degree $2$.
Furthermore they
intersect on $R$ in a curve of type $(2h-f)(3h-2f)=6h^2-7hf$, i.e. a curve
of degree
$11$.  The union of the $5$ singular lines in this intersection is a curve
of type
$5h^2-10hf$, so the remaining part is of type  $h^2+3hf$. This is a curve
of degree
$6$. Again the surface $W_3\cap W_4$ has degree $2$.  If it is irreducible,
then the
lines of one of the rulings are at least $3$-secants to $R$, impossible, so the
quadric surface is contained in $R$, again contradicting our assumption.
If the
surface $W_3\cap W_4$ is the union of two planes, then at least one of the
planes
intersects $R$ in a curve of degree at least $3$, again contradicting our
assumption,
so this case is also impossible.  This concludes the proof that $W$ is
irreducible.
\skipaline
It remains to show that $R$ is contracted by $\varphi$.  If not, $W$ must be
contracted. First of all this is possible only if $W$ is swept out by
lines, i.e.
has a line through every point.  For this consider a general plane in $R$.  It
intersects $W$ in a plane quintic curve with 6 double points. $W$ has
degree $7$,
so in a general $\Pn 3$ through the plane $W$ has a curve of degree $2$
residual
to the plane curve.  Since the arithmetic genus of $W$ is 1 and the plane
curve is
rational, this residual curve is a conic or two connected lines.  Now, $W$ can
only have a line through each point if this residual curve is always two
lines.
On the other hand if these two lines are contracted by $\varphi$, they are
contracted to the same point, since the matrix $M$ has rank $3$ even in the
point
of intersection.  But all fibers of $\varphi$ are linear, so $W$ must contain a
plane through the two lines, this is absurd.  Thus $W$ is mapped
birationally to
$\Pn 3$ by $\varphi$.

\par
 It follows that the restriction of $\varphi$ to $R$ is a contraction onto
a surface.
Clearly it must contract the lines transverse to the planes.  From the
matrix $M$ we see
that the map is defined by cubics on each plane, and it is birational, so
the image has
degree at least $3$.  Therefore, the inverse map is defined by hyperplanes or
quadrics.  Since $W$ has degree $7$ the inverse map must be defined by
quadrics, and
in fact with one basepoint, i.e. $W$ is a Del Pezzo $3$-fold
$W_7$.\qed

\Section 5.  Non-normal Calabi-Yau $3$-folds in $\Pn 5$

We shall show that $W_t$ for $t=5,6,7$ is bilinked on a quintic
hypersurface to a non-normal Calabi-Yau $3$-fold. $W_t$ is assumed to be general
in the sense of section 4,
 including in the case of $W_7$ any Del Pezzo $3$-fold characterized as in Proposition 4.10.  Most of the argument
goes through also for
$W_8$, but since we shall not need this later we do not conclude in this case (see Remark 5.3).

\Lemma 5.1. The general Del Pezzo $3$-fold
$W_t$ for
$t=5,6,7,8$ is contained in an irreducible quintic hypersurface. The
general such
quintic hypersurface is normal, it has double points
along the singular lines of $W_t$ and
has only canonical
singularities.

\Proof  We start with a general $W_t$ and a general quintic hypersurface
$Q$ through $W_t$.  Since $W_t$ has non-normal double points along the singular lines and quintics generate
its ideal, $Q$ has multiplicity $2$ at a general point on a singular line.
  Let
$$p^{\prime}:{\rm Bl}_s(\Pn 5)\to \Pn 5$$
be the blowup of $\Pn 5$ along the singular lines of
$W_t$, and denote by $W^{\prime}_t$ and $Q^{\prime}$ the strict
transforms of $W_t$ and $Q$.   Then
$W^{\prime}_t$ has at most ordinary double points as singularities. \par
 Next we consider the blowup
 $$p^{\prime\prime}:{\rm Bl}_{s,W} (\Pn 5)\to{\rm Bl}_s(\Pn 5)$$
 of ${\rm Bl}_s(\Pn
5)$ along $W^{\prime}_t$.  Over the singular points of
 $W^{\prime}_t$, the
blowup ${\rm Bl}_{s,W} (\Pn 5)$ will have isolated double points.  Denote by
 $Q^{\prime\prime}$ the strict
transform of $Q^{\prime}$  and let $W^{\prime\prime}_t$ be the strict
transform of
$W^{\prime}_t$ in $Q^{\prime\prime}$.  By the fundamental property of
blowup, $W^{\prime\prime}_t$
is a Cartier divisor on $Q^{\prime\prime}$.\par
 We now analyse the
situation, assuming that we have chosen $Q$ general. For $t\leq 7$, the
ideal of
$W_t$ is generated by quintics by Proposition 4.5, so in this case
$Q^{\prime\prime}$
is smooth outside $W^{\prime\prime}_t$. We claim that $Q^{\prime\prime}$
is smooth and that $p^{\prime\prime}$
restricted to $Q^{\prime\prime}$ is a small resolution of $Q^{\prime}$.
Since the quintics
generate the ideal of $W_t$, the conormal sheaf of  $W^{\prime}_t$ in ${\rm
Bl}_{s} (\Pn
5)$ twisted by the class of
$Q^{\prime}$ is generated by its global sections.  In particular, any divisor
equivalent to $Q^{\prime}$ induces a section of this sheaf, which is a rank $2$
bundle outside the singular points. Its zero locus is the subscheme of
$W^{\prime}_t$
defined by the Jacobian ideal of $Q^{\prime}$. For general $Q$ this is a smooth
curve which does not pass through the singular points of $W^{\prime}_t$.
Thus $Q^{\prime}$ is
singular only along some smooth curve on $W^{\prime}_t$ not passing through
the singularities of $W^{\prime}_t$. The map $p^{\prime\prime}$ defines
an isomorphism between $W^{\prime\prime}_t$ and $W^{\prime}_t$. Moreover
since $W^{\prime\prime}_t$ is a Cartier divisor on $Q^{\prime\prime}$ it
follows that $Q^{\prime\prime}$ is smooth and that $p^{\prime\prime}$
defines a small resolution of $Q^{\prime}$ showing that $Q^{\prime}$ has
only canonical singularities.
For $W_8$ one may
show with [MAC] in the example of Lemma 4.9 that the baselocus of the
quintics through $W_8$ is
precisely the singular lines and the two planes.  Therefore also in this
case  $Q^{\prime\prime}$
is smooth and we can argue as before.\par
Finally
$Q$ is smooth in codimension 1
and since it is a hypersurface, it is normal. It remains to prove that
$Q$ has canonical singularities along the singular lines. First recall
that $Q$ has multiplicity $2$ at the generic point of the singular lines.
Moreover the blowup $p^{\prime}$ defines a resolution of $Q$
over the points on the singular lines outside the singular curve
of $Q^{\prime}$. Hence $Q$ has transversal ordinary double points
outside the finitely many points on the singular lines which are the
intersection with the image  under $p^{\prime}$ of the singular curve of
$Q^{\prime}$. In particular $Q$ has canonical singularities outside these
finitely many points.
Combining this with the argument that $p^{\prime\prime}$ defines a small
resolution of $Q^{\prime}$ gives the claim of the lemma.
\qed

\Proposition  5.2.  The general Del Pezzo $3$-fold $W_t$ for $t=5,6,7$ is
bilinked $(5,4)$
and $(5,5)$ on a quintic hypersurface to a variety $Y$ of degree $t+5$
which has
non-normal
singularities along the singular lines of $W_t$.  The normalization of $Y$ is
a smooth Calabi-Yau $3$-fold.

\Proof  We start with a general quintic hypersurface $Q$ which contains
$W_t$ as in Lemma 5.1, and
keep the same setup and notation as in the above proof.  In addition we let
$H_{Q^{\prime}}$ and
$H_{Q^{\prime\prime}}$ be the pullback by $p^{\prime}$ and
$p^{\prime\prime}$ of a hyperplane $H$
restricted to $Q$.   A $3$-fold $Y$ bilinked in hypersurfaces of degree
$4$ and $5$ to $W_t$ on $Q$ is nothing but a Weil divisor equivalent to
$W_t+H_Q$,
where $H_Q$ is the restriction of $H$ to $Q$.  The strict transform
$Y^{\prime\prime}$ of $Y$ on $Q^{\prime\prime}$ is a Cartier divisor
linearly equivalent to
$W_t^{\prime\prime}+H_{Q^{\prime\prime}}$.  To analyze the singularities
and the canonical sheaf of
$Y$ we perform adjunction on the smooth $4$-fold $Q^{\prime\prime}$.\par
Let  $E^{\prime\prime}$ be
the pullback of the exceptional divisor of ${\rm Bl}_{s} (\Pn
5)$ to ${\rm Bl}_{s,W} (\Pn 5)$. Then the canonical line bundle on
$Q^{\prime\prime}$ is ${\cal
O}_{Q^{\prime\prime}}(E^{\prime\prime}-H_{Q^{\prime\prime}})$.
By adjunction on $Q^{\prime\prime}$,
the canonical line bundle on $W_t^{\prime\prime}$ is
${\cal O}_{W_t^{\prime\prime}}(W_t^{\prime\prime}+
E^{\prime\prime}-H_{Q^{\prime\prime}})$.
Consider now the normalization $\tilde W_t$ of $W_t$. The map
$W_t^{\prime\prime}\to W_t$ factors through this normalization. On $\tilde W_t$
there are no exceptional divisors since the singular locus has codimension 2.
Therefore, the map $W_t^{\prime\prime}\to \tilde W_t$ is the blowup along a
curve with
exceptional line bundle ${\cal
O}_{W_t^{\prime\prime}}(E^{\prime\prime})$. But $\tilde W_t$ is a Del Pezzo
$3$-fold with only
canonical Gorenstein
singularities,  so the canonical divisor of $W_t^{\prime\prime}$
is also
${\cal O}_{W_t^{\prime\prime}}(E^{\prime\prime}-2H_{Q^{\prime\prime}}).$
Therefore ${\cal
O}_{W_t^{\prime\prime}}(W_t^{\prime\prime}-H_{Q^{\prime\prime}})={\cal
O}_{W_t^{\prime\prime}}(-2H_{Q^{\prime\prime}})$
and we obtain the equality
$$
{\cal O}_{W_t^{\prime\prime}}(W_t^{\prime\prime}+H_{Q^{\prime\prime}})=
{\cal O}_{W_t^{\prime\prime}}.
$$
Now, consider the exact sequences of sheaves on $Q^{\prime\prime}$
$$
0\longrightarrow {\cal
O}_{Q^{\prime\prime}}(H_{Q^{\prime\prime}})\longrightarrow {\cal
O}_{Q^{\prime\prime}}(W_t^{\prime\prime}+H_{Q^{\prime\prime}})\longrightarrow
{\cal O}_{W_t^{\prime\prime}}(W_t^{\prime\prime}+
H_{Q^{\prime\prime}})\longrightarrow
0,
$$

$$0\longrightarrow {\cal
O}_{Q^{\prime\prime}}\longrightarrow {\cal
O}_{Q^{\prime\prime}}(H_{Q^{\prime\prime}})\longrightarrow {\cal
O}_{H_{Q^{\prime\prime}}}(H_{Q^{\prime\prime}})\longrightarrow 0$$
and
$$0\longrightarrow {\cal
O}_{H_{Q^{\prime\prime}}}\longrightarrow {\cal
O}_{H_{Q^{\prime\prime}}}(H_{Q^{\prime\prime}})\longrightarrow {\cal
O}_{H_{Q^{\prime\prime}}\cap
{H^{\prime}}_{Q^{\prime\prime}}}(H_{Q^{\prime\prime}})\longrightarrow 0$$
for general
hyperplanes $H$ and $H^{\prime}$.  Since $Q$ and $Q\cap H$ have only canonical
singularities, $Q^{\prime\prime}$ and $H_{Q^{\prime\prime}}$ are regular.
Thus the
second and third sequence remain exact on global sections.  Furthermore $H\cap
H^{\prime}\cap Q$ is a smooth quintic surface, so
${\rm h}^1({\cal
O}_{H\cap
H^{\prime}\cap Q}(H))={\rm h}^1({\cal
O}_{H_{Q^{\prime\prime}}\cap
{H^{\prime}}_{Q^{\prime\prime}}}(H_{Q^{\prime\prime}}))=0$.  Therefore
${\rm h}^1({\cal
O}_{H_{Q^{\prime\prime}}}(H_{Q^{\prime\prime}}))={\rm h}^1({\cal
O}_{Q^{\prime\prime}}(H_{Q^{\prime\prime}}))=0$ and also the first sequence
remains exact after taking global sections. Thus the line bundle in the
middle is
generated by global sections, i.e. the linear system
$|W_t^{\prime\prime}+H_{Q^{\prime\prime}}|$ of divisors on
$Q^{\prime\prime}$ has no
basepoints.  By Bertini,
we may conclude that a general member $Y^{\prime\prime}$ of this linear
system is smooth. \par
The exceptional locus of $p^{\prime\prime}$ on $Q^{\prime\prime}$ is a
surface scroll, and
$W_t^{\prime\prime}$ intersects this scroll in a section over the smooth
base curve.   The
restriction of the linear system
$|W_t^{\prime\prime}+H_{Q^{\prime\prime}}|$ to this scroll has no
basepoints so the general member is again a section over the base curve.
Therefore the restriction
of $p^{\prime\prime}$ to $Y^{\prime\prime}$ is also an isomorphism, and the
image
$Y^{\prime}$ on $Q^{\prime}$ is smooth, while $Y$ is singular only along
the singular lines.  Since $W_t$ has non-normal double points along the
singular lines, the
same is the case for $Y$.\par The canonical line bundle on
$Y^{\prime\prime}$ is easily computed by
adjunction on $Q^{\prime\prime}$. In fact it is
$${\cal
O}_{Y^{\prime\prime}}(Y^{\prime\prime}+E^{\prime\prime}-H_{Q^{\prime\prime}})=
{\cal
O}_{Y^{\prime\prime}}(W_t^{\prime\prime}+E^{\prime\prime}).$$
Since ${\cal
O}_{W_t^{\prime\prime}}(Y^{\prime\prime})={\cal
O}_{W_t^{\prime\prime}}(W_t^{\prime\prime}+H_{Q^{\prime\prime}})={\cal
O}_{W_t^{\prime\prime}}$, we have ${\cal
O}_{Y^{\prime\prime}}(W_t^{\prime\prime})={\cal
O}_{Y^{\prime\prime}}$.  Therefore the canonical line bundle on
$Y^{\prime\prime}$ is
${\cal
O}_{Y^{\prime\prime}}(E^{\prime\prime})$.\par
 We turn to the singular variety $Y$. It is singular precisely along the
singular lines of $W_t$ where it has non-normal double points.  Consider
the normalization $\tilde
Y\to Y$.  As in the case of $W_t$, the resolution of singularities
$Y^{\prime\prime}\to Y$ factors
through this normalization.  We will show that $\tilde Y$ is already
smooth.  For this we first
specialize the quintic hypersurface $Q$ to a hypersurface containing $W_t$
and a $3$-fold $Z$ that
is smooth along the singular lines of $W_t$. In case $t=7$, we take $Z$ to
be the Segre cubic
scroll $R$, while in the cases $t=5$ and $t=6$ we may take $Z$ to be one
and two $\Pn 3$'s
respectively. We make the computation explicit in the case when
$t=7$, which we will use later. The other cases are similar. Assume
that $Q$ is a quintic
which contains $W_t$ and the rational cubic scroll $R$.  Let $Q_0\to Q$ be the
blow up of $Q$ along $R$.  This defines a small resolution of the
singularities of $Q$ along the
lines. Over each line there is an exceptional scroll.  This scroll is
isomorphic to some Hirzebruch
surface $F_a$ for some $a\geq 0$.  The strict transform of $W_t$ meets this
scroll in a rational
curve which is a bisection on  the scroll, since $W_t$ is double along the
line.  Since $Y$ as a
Weil divisor is equivalent to $W_t+H_Q$ on $Q$, and the strict transform of
$H_Q$ intersects the
scroll in a ruling, we get that the strict transform of $Y$ meets the
exceptional scroll in a
bisection which is an elliptic curve. The normalization of $Y$ factors
through this small
resolution and, therefore, has elliptic curves lying over the singular
lines.   In particular the
normalization is smooth over the singular lines.  By deformation to the
general quintic $Q$, the
normalization of $Y$ is smooth.\par Now, the normalization $\tilde Y$ of
$Y$ is isomorphic to $Y$
in codimension 1, therefore the canonical line bundle ${\cal
O}_{Y^{\prime\prime}}(E^{\prime\prime})$ on $Y^{\prime\prime}$ is the
exceptional line bundle of the
map $Y^{\prime\prime}\to \tilde Y$ and the canonical line bundle of $\tilde
Y$ is trivial.\par The
irregularity of $\tilde Y$ equals the irregularity of $Y^{\prime\prime}$,
since they are birational
and both are smooth.   But on $Y^{\prime\prime}$ we get $${\rm h}^1({\cal
O}_{Y^{\prime\prime}})=\coh
2{Q^{\prime\prime}}{-W_t^{\prime\prime}-H_{Q^{\prime\prime}}}=\coh
1{W_t^{\prime\prime}}{-H_{Q^{\prime\prime}}}=0$$
 so the normalization of $Y$ is regular.  Thus the
normalization of $Y$ is a smooth Calabi-Yau $3$-fold and the proposition
follows. \qed

\Remark 5.3. The above proof goes through also for $W_8$, except
for the question of a smooth normalization. \par

\Section 6. Equations of elliptic scrolls in $\Pn 5$

In the remaining sections we prove that the union of two elliptic scrolls which
intersect along
an abelian surface is bilinked to a Del Pezzo $3$-fold $W_7$.  On the way we
describe the
family of these reducible Calabi-Yau $3$-folds using Heisenberg symmetry.
But first we
study the ideal of elliptic scrolls without using this symmetry.
   \medskip Recall from section 3 that two elliptic
scrolls whose intersection is a $(1,6)$-polarized abelian surface in $\Pn
5$, are each
singular along three lines.  Therefore we restrict our attention to this
kind of elliptic
scrolls.  We start with a lemma which gives a quick
construction of such scrolls.

\Proposition 6.1. The union of three $\Pn 3$'s in
$\Pn 5$, which meet pairwise in lines, is linked $(3,3)$ to an elliptic
scroll.  Furthermore any elliptic $3$-fold scroll
 of degree $6$, singular along
three lines which span $\Pn 5$, is linked $(3,3)$ to three
$\Pn 3$'s.\par  \Proof
First, consider three $\Pn 3$'s which meet pairwise in lines and two general
cubic hypersurfaces containing them.  The linked variety $X$ has degree $6$ and
sectional genus 1 (cf. [PS]).  Since the complete intersection has trivial
canonical bundle, each component intersects the rest along an anticanonical
divisor.  Therefore $X$ meets each of the $\Pn 3$'s in a quartic surface,
singular along the two lines of intersection with the other two $\Pn 3$'s.
These quartic surfaces clearly are elliptic scrolls:
Through every point on the surface outside the two
singular lines there is a unique line in the surface
through the point intersecting the two lines, by
Bezout, so the quartic surface is a scroll of lines.  Through every
point on the singular lines there are two rulings of
the scroll.  In the Grassmannian of lines the curve
parametrizing the rulings has a double cover to each
of the two singular lines, so the curve is of type
$(2,2)$ on $\Pn 1\times \Pn 1$, i.e. it is elliptic.
The normalization of $X$ therefore contains elliptic
scrolls in parts of hyperplane sections.  On the other
hand, residual to each quartic surface in hyperplane
sections, there is a pencil of surfaces of degree two on
$X$. Residual to the $\Pn 3$ they form part of a complete
intersection of two quadrics in the hyperplane.  The
other two $\Pn 3$'s intersect this hyperplane in two
planes passing through a point, so the surface of degree
two on $X$ must be two planes residual to these in the
complete intersection. This displays the scroll structure
of $X$.  Clearly the scroll is elliptic. \par Next,
consider an elliptic scroll $X_F$ of degree $6$ and
singular along three lines that span $\Pn 5$. Let
$$\varphi : {\bf P}({\cal E})\to X_F$$ be the
normalization map. Like in section 3 the three lines
correspond to a decomposition
 $${\cal
E}={\cal L}_0\oplus {\cal L}_1\oplus {\cal L}_2,$$ where the
${\cal L}_i$  are line bundles of degree $2$. We have
$${\rm H}^0(F,{\cal E})\cong {\rm H}^0(X_F,{\cal O}_{X_F}(1))\cong {\rm
H}^0(\Pn 5,{\cal O}_{\Pn 5}(1))$$
since the linear system that maps ${\bf P}({\cal E})$ into $\Pn 5$ is complete.
There is an isomorphism
$${\rm H}^0({\cal O}_{{\bf P}({\cal E})}(n))\cong {\rm H}^0(F,{\rm S}^n{\cal
E})\cong
{\rm H}^0(F,{\cal L}_0^n\oplus {\cal L}_0^{n-1}\otimes {\cal
L}_1\oplus\ldots\oplus {\cal L}_2^n)$$
where the number of summands is
$${{n+2}\choose 2}={{(n+1)(n+2)}\over 2}$$
and the degree of each summand is $2n$.
This shows that
$$h^0({\cal O}_{{\bf P}({\cal E})}(n))=h^0(F,{\rm S}^n{\cal E})=n(n+1)(n+2).$$
We now compare this to the situation on the singular scroll $X_F$.  Recall
that there
are three sections $F_i\subset {\bf P}({\cal E})$ such that the map
$\varphi$  restricted
to $F_i$ induces a double cover
$\varphi :F_i\to L_i$ onto a line $L_i$.  This is branched over 4 points
and hence
$$\varphi_{\ast}{{\cal O}_{F_i}}={\cal O}_{L_i}\oplus {\cal O}_{L_i}(-2),$$
and we have the exact sequence
$$0\to {\cal O}_{L_i}\to \varphi_{\ast}\varphi^{\ast}{\cal O}_{L_i}\to {\cal
O}_{L_i}(-2)\to 0$$
where $\varphi_{\ast}\varphi^{\ast}{\cal O}_{L_i}=\varphi_{\ast}{\cal
O}_{F_i}$.

\Lemma 6.2.  There is an exact sequence of sheaves
$$0\to {\cal O}_{X_F}\to \varphi_{\ast}{\cal O}_{{\bf P}({\cal E})}\to
\oplus_{i=0}^2{\cal
O}_{L_i}(-2)\to 0.$$
\par
\Proof  The cokernel clearly has support on the locus on $X_F$ where the
map is not
an isomorphism, i.e. on the three lines $L_i$.  Restricted to these lines the
exact sequence reduces to the one above.  Since the lines are disjoint, the
cokernel is a
direct sum.\qed

\Lemma 6.3.  After tensoring with ${\cal O}_{X_F}(n)$ for any nonnegative
$n$ the
sequence of Lemma 6.2 is exact on global sections.

\Proof For $n=0,1$ the statement is immediate, since the third term has
no sections.
Furthermore  $\coh 1{X_F}1={\rm h}^1({\cal O}_{{\bf P}({\cal E})}(1))=0$ and
$\coh 2{X_F}n={\rm h}^2({\cal O}_{{\bf P}({\cal E})}(n))=0$ for $n\geq 1$.
Let $H$ be a general hyperplane, and $P=\Pn 3$ be a general $3$-space inside
$H$. Then ${X_F}\cap P$ is a smooth elliptic curve, so $\coh 1{{X_F}\cap
P}n=0$ for $n\geq 1$.  Furthermore
$\coh 1{{X_F}\cap H}1=0$.  Inductively $\coh 1{{X_F}\cap
H}n\leq\coh 1{{X_F}\cap H}{n-1}=0$ for $n\geq 2$, and similarly
$$\coh 1{{X_F}}n\leq\coh 1{{X_F}}{n-1}=0$$
for $n\geq 2$ and the lemma follows.\qed

 Lemma 6.3 allows us to compute
$$h^0({\cal O}_{X_F}(n))=n(n+1)(n+2)-3(n-1), \quad n\geq 1.$$
The natural maps
${\rm H}^0({\cal O}_{\Pn 5}(n))\to {\rm H}^0({\cal O}_{X_F}(n))$
and the dimension
$$h^0({\cal O}_{\Pn 5}(n))={\rm dim}\ {\rm S}^n{\rm
H}^0({\cal E})={{n+5}\choose 5},$$
gives us
$$\cohi 0{X_F}n-\cohi 1{X_F}n ={{n+5}\choose 5}+3(n-1)-n(n+1)(n+2).$$
In particular $\cohi 0{X_F}3\geq 2$.  Clearly the complete intersection of two
cubics through $X_F$ contains the $\Pn 3$'s spanned by the pairs of singular
lines, therefore \break
$\cohi 0{X_F}3=2$ and Proposition 6.1 follows.\qed

\Corollary 6.4.  $\cohi 0{X_F}5=54$.\par

\Proof  $$\cohi 0{X_F}5-\cohi 1{X_F}5 =54,$$
as follows from the computation in the proof of Proposition 6.1.  But
$X_F\cap \Pn 2$ for a general $\Pn 2$ is $4$-regular in the sense of
Castelnuovo-Mumford, so $\cohi 1{X_F}n=0$ for $n\geq 3$.\qed

\Section 7. Heisenberg symmetry of elliptic scrolls in $\Pn 5$

First we collect some basic observations.
Recall the Heisenberg group of level 6:
$$
1\rightarrow \mu_6\rightarrow H_6\rightarrow\Z/6\times
\Z/6\rightarrow 0.
$$
Given any subgroup $G\subset \Z/6\times\Z/6$ we can consider the
preimage
$$
H(G):=\pi^{-1}(G)\subset H_6.
$$\
Let $\langle x_0,\ldots,x_5\rangle$ be a basis for a $6$-dimensional
complex vector space $V$ and let $\langle e_0,\ldots,e_5\rangle$ be a basis
for $V^*$. The Schr\"odinger representation $\rho:H\rightarrow {GL}(V,\C)$ is
defined by $$\sigma: x_i\mapsto x_{i+1},\quad \tau: x_i\mapsto \rho^{i}x_i$$
where $\rho =\rho=e^{2\pi i/6}$ and indices are taken modulo $6$.  It defines,
by restriction, representations $$ \rho_G:H(G)\rightarrow {GL}(V,\C). $$
The group $\Z/6 \times\Z/6$ contains unique subgroups
isomorphic to  $\Z/2\times \Z/2$, resp.
$\Z/3\times\Z/3$. The preimages of these groups in
$H_6$ are isomorphic to the Heisenberg groups $H_2$ and $H_3$ of level $2$ and
$3$. Note, however, that in the case of $H_3$ the induced representation
differs
from the Schr\"odinger respresentation by the non-trivial automorphism of the
Galois group.

\Lemma 7.1. {\rm(i)} There are $4$ subgroups isomorphic to $\Z/3$ in
$\Z/3\times\Z/3$.
\smallskip {\rm(ii)} There are $4$ subgroups isomorphic to
$\Z/2\times\Z/6$ in $\Z/6\times\Z/6$.

\Proof  Claim {\rm(i)} is trivial.
Every subgroup isomorphic to
$\Z/2\times\Z/6$ is generated by the elements of order $2$ and one
element of order $3$ in $\Z/6\times\Z/6$, i.e. by the subgroup
$\Z/2\times \Z/2\subset \Z/6\times\Z/6$ and a subgroup of
order $3$ of the group
$\Z/3\times\Z/3\subset\Z/6\times\Z/6$. This shows {\rm(ii)}.
\hfill\qed
\skipaline

\Remark 7.2.  At the same time this gives us a natural $1:1$ correspondence
between the $4$ groups of part {\rm(i)} of Lemma 7.1 and the 4 groups of part
{\rm(ii)}. We shall use this frequently in
what follows.

We shall denote the 4 subgroups of $\Z/6\times \Z/6$ which are isomorphic to
$\Z/2\times\Z/6$ by $K_1,\ldots, K_4$.
For every group $G$ in $\Z/6\times \Z/6$ we set
$$
H(G)_{\iota} = \langle H(G), \iota\rangle \subset G_6.
$$
We denote by $G_2$ resp. $G_3$ the groups $\langle H_2,\iota\rangle$, resp.
$G_3=\langle H_3,\iota\rangle.$ Next we need some elementary representation
theory. We denote by $U$ the Schr\"odinger representation  of
$H_2$.

\Lemma 7.3. {\rm(i)} As an $H_2$-module $V\cong 3U$,\smallskip
{\rm(ii)} as a $G_2$-module $V\cong 2 U_+\oplus U_-$, where $U_+$, resp. $U_-$
means that $\iota$ acts by $+1$, resp. $-1$ on $U$,\smallskip
{\rm(iii)} as $K_i$-modules $V\cong U\oplus U'\oplus {\bar {U}}'$. Here the
subgroup of $K_i$ which is isomorphic to $\Z/3$, acts on $U'$ by a non-trivial
character and on ${\bar U}'$, by the inverse of this character. The involution
$\iota$ leaves $U$ fixed and interchanges $U'$ and ${\bar{U}}'$.

\Proof (i), (ii) We can decompose $V$ as an $H_2$-module as
follows:
$$
\langle x_0, x_3 \rangle,\langle x_2, x_5 \rangle, \langle x_4, x_1 \rangle.
$$
From this the claim is obvious.\smallskip
(iii) It is enough to consider one of the groups $K_i$. The others can be done
in the same way, or one can use the normalizer $N_6$ of $H_6$ in $GL(6,\C)$.
Here we shall consider the subgroup $K$ given by $\tau^2$. One
immediately checks that $\tau^2$ acts by $1,\rho^4,\rho^2 (\rho=e^{2\pi
i/6})$ on the above submodules of $V$, hence giving the claim.
\hfill\qed

We next want to associate basic geometric objects in $\Pn 5$ to the subgroups
$K_i, i=1, 2, 3, 4$.

\Lemma 7.4. {\rm(i)} To every subgroup $H(K_i), i=1, 2, 3, 4$ of $H_6$
one can associate a unique set of $3$ lines in $\Pn 5$ which is an $H_6$-orbit
such that $H(K_i)$ is the stabilizer of each of these lines.  The
distinguished subgroup of
order $3$ in $K_i$ fixes these lines pointwise. \smallskip {\rm(ii)} Every
$H_6$-orbit of lines
in $\Pn 5$ consisting of $3$ lines is one of the above.

\Proof Let $\{L_1, L_2, L_3\}$ be an $H_6$-orbit of lines in $\Pn 5$. Then
every line $L_j$ has a stabilizer in $\Z/6\times \Z/6$ of order $12$. This must
then be one of the groups $K_i, i=1, 2, 3, 4$.
Without loss of generality it suffices to consider the
group generated by $\langle\tau,\sigma^3\rangle$. Then the
action of $\tau$ on the vector space associated to such a
line $L_j$ splits into a sum of two different characters.
Therefore $L_j$ is spanned by two basis vectors $e_k,
e_l$. To obtain invariance under $\sigma^3$ the only
possibilities are $\langle e_0, e_3\rangle, \langle e_1,
e_4\rangle$ and $\langle e_2, e_5\rangle$. Furthermore the
distinguished subgroup order 3 generated by $\tau^2$
fixes these three lines pointwise. \hfill\qed

\Lemma 7.5. {\rm(i)} To every subgroup $H(K_i), i=1,2,3,4$ of $H_6$ one
can associate a unique set of three $3$-spaces in $\Pn 5$ which is an
$H_6$-orbit
such that $H(K_i)$ is the stabilizer of each of these ${\Pn 3}'{\it
s}$.\smallskip {\rm(ii)} Every $H_6$-orbit of $3$-spaces in $\Pn 5$ consisting
of three ${\Pn 3}'{\it s}$ is one of the above.

\Proof This is the dual statement to Lemma 7.4. Given three
lines $L_1, L_2, L_3$
which form an $H_6$-orbit, the three ${\Pn 3}'{\rm s}$ are the spaces
spanned by
two of these lines.
\hfill\qed

Next we turn to the space of cubic forms $\Coh 0{\Pn 5}3$ which
we want to study as an $H_6$-, resp. $G_6$-module.
\Lemma 7.6. The $G_6$-module $\Coh 0{\Pn 5}3$ is a sum of four
$2$-dimensional and twelve $4$-dimensional representations.  As
$H_3$-representation it is a sum of characters. The
trivial character corresponds to the four pencils. The other $8$ come in
pairs (given
by the involution $\iota$) and each pair determines three $4$-dimensional
irreducible
$G_6$-representations. The subspace of
$2$-dimensional representations is spanned by
$$\matrix {\langle x_0^3+x_2^3+x_4^3, x_1^3+x_3^3+x_5^3\rangle\cr
\langle x_0x_2x_4, x_1x_3x_5\rangle\cr
\langle x_3^2x_0+x_5^2x_2+x_1^3x_4, x_4^2x_1+x_0^2x_3+x_2^3x_5\rangle\cr
\langle x_1x_2x_3+x_3x_4x_5+x_5x_0x_1,
x_2x_3x_4+x_4x_5x_0+x_0x_1x_2\rangle.\cr}
$$

\Proof  This is a straightforward computation.\qed

Since all the $2$-dimensional representations are mutually isomorphic this
defines
a $\Pn 3$ of
pencils of cubics.

\Proposition 7.7.
For every $G_6$-orbit of 3 lines in $\Pn 5$, there is a
unique pencil of $G_6$-invariant pencils of cubics containing these lines.

\Proof We can assume that the 3 lines in question are $\langle e_0, e_3
\rangle,
\langle e_1, e_4 \rangle$ and $\langle e_2, e_5 \rangle$. A general pencil of
$G_6$-invariant pencils of cubics is of the form

$$\matrix
{a\langle x_0^3+x_2^3+x_4^3,x_1^3+x_3^3+x_5^3\rangle\cr
+b\langle x_0x_2x_4,x_1x_3x_5\rangle\cr
+c\langle x_3^2x_0+x_5^2x_2+x_1^2x_4,x_4^2x_1+x_0^2x_3+x_2^2x_5\rangle\cr
+d\langle
x_1x_2x_3+x_3x_4x_5+x_5x_0x_1,x_2x_3x_4+x_4x_5x_0+x_0x_1x_2\rangle
.\cr}$$ Such a pencil contains the above lines if and only if
$a=c=0$. \hfill\qed

\Remark 7.8. {\rm(i)} It is also easy to determine the pencil of pencils
containing the other minimal $H_6$-orbits of lines. E.g. the three lines fixed
by the elements $\sigma$ and $\tau^3$ are the line
$\langle (1,1,1,1,1,1), (1,-1,1,-1,1,-1) \rangle$ and its
$\tau$-translates. The corresponding pencil of pencils is given by
$3a+b=c+d=0$.\par
{\rm(ii)}
Every pencil of $G_6$-invariant pencils of cubics containing a minimal orbit
$\{L_1, L_2, L_3\}$ also contains the three ${\Pn 3}'s$ spanned by two of these
lines. (This can be seen by direct inspection). Hence every such pencil has a
base locus consisting of three ${\Pn 3}'s$ and a residual $3$-fold $X$ of
degree $6$.   By Proposition 6.1, $X$ is an elliptic
scroll.  Furthermore, every elliptic $3$-fold scroll $X_F$, as in section 3,
is $G_6$-invariant, so it is contained in a
$G_6$-invariant pencil of cubics. Hence for a general
pencil of cubics the residual $X$ must be of the form
$X_E$ for a suitable elliptic curve $E$.
\skipaline
Our next aim is to study the Heisenberg action on the embedded abelian surfaces
$E \times F$. Recall that every abelian surface with a very
ample $(1,6)$-polarization can be embedded $G_6$-equivariantly into ${\Pn 5}$
and that the choice of such an embedding is equivalent to the choice of a
level-$6$ structure on $E \times F$ (i.e. a canonical level struture
associated to the polarization $H$ by which we mean a
symplectic basis of the kernel of the map $\lambda_H: A
\to \hat A$). We denote the family of all $G_6$-invariant
abelian surfaces in ${\Pn 5}$ with two plane elliptic fibrations by ${\cal A}$.
As in section $1$ we start with curves $E$ and $F$ and a $3:1$ morphism
$\gamma:E\rightarrow F$. We can assume that
$E=\C/(\Z\tau+\Z)$, $F=\C/(\Z  3\tau+\Z)$ and that the map $\gamma$
is induced by $\gamma(z)=3z$. We denote the
generators $\tau$ and
$1$ of the lattice $\Z\tau+\Z$ by $e_1$ and $e_3$ and set $s_6=e_1/6$
and $t_6=e_3/6$. Moreover we denote the generators $3 \tau$ and $1$ of the
lattice $\Z3\tau+\Z$ by $e_2$ and $e_4$.
Then
$\gamma(e_1)=e_2, \gamma(e_3)=3e_4$. The point $u_6=\gamma(s_6)$
is represented by
$e_2/6$ and $\gamma(t_6)=3 v_6$ is represented by $e_4/2$. We choose $v_6$ as
the point represented by $e_4/6$. In this set-up the product $A=E\times F$ is
given by the period matrix
$$
\left(
\matrix{
\tau & 0     & 1 & 0\cr
0    & 3\tau & 0 & 1\cr}
\right).
$$
The first step is to understand the polarization $H=2E+\Gamma$ in terms of a
Riemann form. This is necessary to understand the level-$6$ structures on $A$.
First we look at the semi-positive line bundle defined by $E$. The
corresponding Riemann form with respect to
 the lattice $L=\Z e_1+\Z e_2+\Z e_3
+\Z e_4$ is clearly given by
$$
H_E=\left(\matrix{
0 & 0 & 0 & 0\cr
0 & 0 & 0 & 1\cr
0 & 0 & 0 & 0\cr
0 & -1 & 0 & 0\cr}
\right).
$$

Next we want to identify the form $H_{\Gamma}$ with respect to the chosen basis
$e_1,\ldots , e_4$. Since $\Gamma\cdot F=1$ we can write $A=\Gamma\times F$.
By abuse of notation, let
$\gamma=(\rm{id},\gamma):E\rightarrow E\times F$ by the
embedding of $E$ into $A$. Then $\gamma(E)=\Gamma$. We
have $\gamma(e_1)=e_1+e_2=:f_1,\gamma(e_3)=e_3+3
e_4=:f_2$. We can also choose $f_1, f_2, e_2, e_4$ as a
basis for $L$. With respect to this basis the
semi-positive form $H_{\Gamma}$ is given by $(e_2,
e_4)=1$ and all other products $0$. A straightforward
calculation then shows that in terms of the basis $e_1,
\ldots, e_4$ the form $H_{\Gamma}$ is given by: $$
H_{\Gamma}=\left(\matrix{ 0 & 0 & 3 & -1\cr
0 & 0 & -3 & 1\cr
-3 & 3 & 0 & 0\cr
1 & -1 & 0 & 0\cr}
\right).
$$

Since $H=2E+\Gamma$ it follows that the corresponding form with respect to the
basis $e_1, \ldots, e_4$ is given by
$$
H=\left(\matrix{
0 & 0 & 3 & -1\cr
0 & 0 & -3 & 3\cr
-3 & 3 & 0 & 0\cr
1 & -3 & 0 & 0\cr}
\right).
$$
Note that this is indeed the form associated to a $(1,6)$-polarization, since
$$
\rm{det}\left(\matrix{3 & -1\cr
-3 & 3\cr}\right) = 6.
$$

Our next aim is to identify the group $\Theta(H)=\ker(\lambda_H:
A\rightarrow{\hat A})$ as a subgroup of $A^{(6)}=E^{(6)}\times F^{(6)}$.
What we have to do is to find a basis of $L^{\vee}/L$ where $L^{\vee}$ is
the dual lattice
with respect to the form $H$. General theory tells us that $L^{\vee}/L\cong
\Z/6\times\Z/6$. It is a straightforward calculation to check that $e_1/2 +
e_2/6, e_3/6 + e_4/2\in L$. We can take these elements as generators of
$L^{\vee}/L$. Note that as points in $A=E\times F$ these are just the points
$(3s_6, u_6)$ and $(t_6, 3v_6)$.

Let $w_1=e_1/2+e_2/6, w_2=e_3/6+e_4/2$. Then a straightforward
calculation shows that for the Weil pairing with respect to $H$ we have
$$
(w_1, w_2) = \left( e^{2\pi i/6}\right),
$$
i.e. these points define a level-$6$ structure. Note that $2 w_1
\in F = \{0\}\times F$ and $2 w_2\in E=E\times\{0\}$. In particular
the $2$ groups of order $3$ generated by $2 w_1$ and $2 w_2$ each respect
one of the
two plane elliptic fibrations of $A = E \times F$. Since the embedded
abelian surface $A=E \times F$ is $G_6$-invariant the same holds for the
scrolls $X_E$ and $X_F$ defined by the $2$ plane cubic fibrations.
Since the groups of order $3$ generated by $2 w_1$, resp. $2 w_2$
each respect one of these fibrations it follows that they act trivially
on the $3$ singular lines of the scroll $X_E$, resp. $X_F$. In
particular this gives $2$ groups $H(K_i)$ and $H(K_j)$ which each has
the $3$ singular lines of the scrolls $X_E$ , resp. $X_F$ as one of its
orbits. Any other choice of a level-$6$ structure on $A$ gives an
analogous picture. We shall return to this in a moment.
\skipaline

We now want to understand the variety ${\cal A}$ parametrizing $G_6$-invariant
abelian surfaces with $2$ plane cubic fibrations. We have already observed
that each of these two fibrations determines a singular scroll $X_E$,
resp. $X_F$ and a group $H(K_i)$, resp. $H(K_j)$. The abelian surface $A$
is the intersection of the scrolls $X_E$ and $X_F$ (cf. Proposition 3.5.) This
defines a decomposition of the family ${\cal A}$ into six families ${\cal
A}_{ij}$ where $\{i,j\}\subset\{1, 2, 3, 4\}$. We want to exhibit a concrete
parametrization of the families ${\cal A}_{ij}$ , thereby also showing that
the ${\cal A}_{ij}$ form six irreducible components. To do this we go back
to an elliptic
curve $E$ as before and the level-$6$ structure on $E$ given by
$(s_6, t_6)$. We can
perform the above construction and associate to these data the surface
$A=E\times F$, the polarization $H=2E+\Gamma$ and the level-$6$ structure
$(w_1, w_2)$. This gives us a morphism
$$
{\psi}_{ij}: X^0(6)\rightarrow{\cal A}_{ij}
$$
from the (open) elliptic modular curve $X^0({6})$ parametrizing
elliptic curves with a level-$6$ strucure to ${\cal A}_{ij}$.
Note that the (compact) modular curve $X(6)$ is an elliptic curve.

\Lemma 7.9. The map ${\psi}_{ij}: X^0(6)\rightarrow{\cal A}_{ij}$ is
surjective onto the component ${\cal A}_{ij}$ and has degree $3$.

\Proof Here we shall treat
the case where $K_i$ and $K_j$ are the groups determined by $\langle \sigma^2
\rangle$ and $\langle \tau^2 \rangle$. This is no loss of generality.
Going back to the abelian surface $A=E\times F$
we want to study the possible embeddings of $A$ into
$\PP^5$ such that $A\in {\cal A}_{ij}$. Since $G_6$-invariant embeddings of
$A$ correspond to the choice of a level-$6$ structure $(w'_1, w'_2)$ we have to
look for those level-$6$ structures $(w'_1, w'_2)$ such that $2w'_1\in
F=\{0\}\times F$ and $2w'_2\in E=E\times\{0\}$ or $2w'_1\in E$ and $2w'_2\in
F$. These two cases correspond to changing the role of $E$ and $F$ and it is,
therefore, enough to look at the first possibility, namely $2w'_1\in F$ and
$2w'_2\in E$. Since $(w_1, w_2)$ is a basis of $\Theta(H)$ we can write
$w'_1=\alpha w_1+\beta w_2, w'_2=\gamma w_1+\delta w_2$. Moreover $2w_1\in F$
and hence $2w'_1\in F$ if and only if $2\beta w_2\in F$ which is only the
case
for $\beta=0$ or $3$. Moreover $w'_1$ has to be an element of order $6$. This
gives us three possibilities for $\pm w'_1$, namely $w_1, w_1+3w_2, 2w_1+3w_2$.
Since $(w'_1, w'_2)$ and $(-w'_1, -w'_2)$ define the same level-$6$ structure
we can assume that $w'_1$ is one of the 3 points above. A similar argument can
be applied to $w'_2$ and altogether we find the following $6$ possibilities
for level-$6$ structure $(w'_1, w'_2)$ with $2w'_1\in F$ and $2w'_2\in E$:
$$
\matrix{
(w_1, w_2),(w_1, w_2+3w_1), (w_1+3w_2, w_2),\cr
(-2w_1+3w_2, w_2+3w_1), (-2w_1+3w_2, -2w_2+3w_1), (w_1+3w_2, -2w_2+3w_1).}
$$
For every pair $(w'_1, w'_2)$ as above the pair $(-w'_2, w'_1)$ is a level-$6$
structure with $-2w'_2 \in E, 2w'_1\in F$. In this way we obtain $12$
level-$6$ structures which belong to the pair $\{i,j\}$. This fits in with the
number of level-$6$ structures which is given by
$$
{1\over2} 6^3 \left(1-{1\over4}\right)\left(1-{1\over 9}\right)=72=12\cdot 6,
$$
where $6$ corresponds to the number of pairs $\{i,j\}\subset\{1, 2, 3, 4\}$.
\skipaline
The map ${\psi}_{ij}: X^0(6)\rightarrow{\cal A}_{ij}$ was defined by
associating to an elliptic curve $E$ with level-$6$ structure $(s_6,
t_6)$ the
abelian surface $A=E\times F$ with level-$6$ structure $(w_1, w_2)$.
We want to prove that
this map is surjective which implies in particular that ${\cal A}_{ij}$ is
irreducible. So we have to show that we can obtain all abelian surfaces
$A=E\times F$ and all level-$6$ structures $(w'_1, w'_2)$ as
above by varying $E$ and the level-$6$ structure $(s_6, t_6)$. Let
$(s'_6, t'_6)$ be any level-$6$ structure on $E$. To this we
associate a surface $A'=E\times F'$ with $F'=E/\langle
2t'_6\rangle$. If we want that $F'=F$ we must (at least for general
$F$) have that $\langle 2t'_6\rangle=\langle 2t_6 \rangle$.
Moreover $t'_6$ must be a point of order $6$. Up to sign this
leaves us with the possibilities $t'_6=t_6, 3s_6+t_6, 3s_6+2t_6$.
Altogether we obtain $18$ possible level-$6$ structures, namely $$
\matrix{(s_6+2it_6, t_6), (s_6+(2i+1)t_6, t_6)\cr (s_6+2it_6,
3s_6+t_6), (-2s_6+(2i+1)t_6, 3s_6+t_6)\cr (-2s_6+(2i+1)t_6,
3s_6-2t_6), (s_6+(2i+1) t_6, 3s_6-2t_6)} $$
where in each case $i=0,1,2$. Each of the values $i=0,1,2$ gives the same
level-$6$ structure on $A$.
Hence under the morphism ${\psi}_{ij}: X^0(6)\rightarrow{\cal
A}_{ij}$ the pairs $(E,(s'_6, t'_6))$ are mapped $3:1$ to
$(A=E\times F, (w'_1, w'_2))$ where $(w'_1, w'_2)$ runs through all
6 possible level-$6$ structures on $A$ with $2w'_1\in F$ and
$2w'_2\in E$. In particular the map $X^0(6)\rightarrow {\cal
A}_{ij}$ is $3:1$ and surjective and ${\cal A}_{ij}$ is
irreducible.
\hfill\qed
\skipaline

Next we consider the family ${\cal V}$ of $G_6$-invariant scrolls which arise
as $\Pn 2$-scrolls defined by a plane cubic fibration of a $G_6$-invariant
abelian surface $A$. Such a scroll is singular along $3$ lines which form an
orbit of one of the groups $H(K_i)$ and hence there is a natural decomposition
${\cal V}= {\cal V}^1 \cup {\cal V}^2 \cup {\cal V}^3 \cup {\cal V}^4$ where
${\cal V}^i$ is the set of those scrolls which are invariant under $H(K_i)$.
These scrolls are in $1:1$ correspondence with points in an open set of the
pencil of pencils of cubics associated to the group $H(K_i)$
(cf. Remark 7.8. ii).
In particular the varieties
${\cal V}^i$ are irreducible and rational.
\medskip

The $G_6$-action on the elliptic scrolls restricts to pencils of
$G_6$-invariant abelian
surfaces. We describe these before we return to the surfaces which lie on
two scrolls.
\medskip
Recall that the normalization $\varphi: {\bf P}_F({\cal E})\to X_F\subset
\Pn 5$ of the scroll
$X_F$ is a $\Pn 2$-bundle over the elliptic curve $F$ associated to the
rank $3$ vector bundle
$${\cal E}={\cal L}_0\oplus h^*{\cal L}_0\oplus (h^2)^*{\cal L}_0,\qquad
(h^3={\rm id}).$$
$X_F$ is singular along 3 lines $L_i$, and
$$(h^i)^*{\cal L}_0=\varphi^{\ast}{\cal O}_{L_i}(1),\qquad i=0,1,2$$ are line
bundles of degree
$2$ on $F$.
If we assume that say
$$L_0=\langle e_0,e_3\rangle,\quad L_1=\langle e_1,e_4\rangle,\quad {\rm
and}\quad L_2=\langle e_2,e_5\rangle,$$ then $${\rm Stab}(L_i)=\langle\tau,
\sigma ^3\rangle.$$ Recall from
section $3$ that we may find ${\bf P}_F({\cal E})$ as a quotient
of $\Pn 2\times E$ by a subgroup ${\Z}/3\subset H_3$.
 This subgroup leaves three
sections $E\to \Pn 2\times E$ invariant, and these three sections are
mapped in the quotient
to the three sections $F\to {\bf P}_F({\cal E})$ which again are mapped $2:1$
to the singular
lines in $X_F$.
In this set-up we may describe the $G_6$-action on the abelian surfaces on
$X_F$.
These surfaces are all pulled back to anticanonical divisors on
${\bf P}_F({\cal E})$.  Recall
that $\coh 0{{\bf P}({\cal E})}{-K}=4.$  Elements in $|-K|$ pull back to
products
$$E^{\prime}\times E\subset \Pn 2\times E$$
where $E^{\prime}$ is a cubic curve invariant under the subgroup of order $3$;
in suitable coordinates it is defined by a form in the web
$$\langle x_0^3, x_1^3, x_2^3, x_0x_1x_2 \rangle.$$
\medskip
These forms are precisely the invariants of degree $3$ of a subgroup of order
$3$ in $H_3$.
In fact the elements of order $3$ of $H_6$ which fix the lines $L_i$
pointwise, leave each
plane in $X_F$ invariant and the lift to an action on
${\bf P}_F({\cal E})$ and $\Pn 2\times E$ which
leaves each plane invariant and fixes the three special
sections pointwise.
These sections, in suitable coordinates, meet each plane
in the points $(1,0,0), (0,1,0), (0,0,1)$ and the action is given by
$\tau\in H_3$ in the
plane.
\skipaline
The involution  $\iota$ leaves $L_0$ fixed while $L_1$ and $L_2$ are
interchanged. This also
lifts to $\Pn 2\times E$.  The  $\iota$-invariant plane cubics in
the above web form the net $$\langle x_0^3, x_1^3+x_2^3, x_0x_1x_2
\rangle.$$ $G_6$ finally permutes the three lines cyclically, so
the corresponding plane cubics are
defined by invariant forms in the variables permuted cyclically.  This
action is
the one defined by $\sigma\in H_3$. The $G_6$-invariant anticanonical
divisors therefore correspond precisely to the Hesse-pencil
 $$\langle x_0^3+x_1^3+x_2^3, x_0x_1x_2 \rangle.$$
\medskip
 As an $H_3$-module $\Coh 0{\Pn 2}3$ splits as the
Hesse-pencil plus 8 characters.  Two of these give rise to anticanonical
divisors on
${\bf P}_F({\cal E})$.
The remaining 6 characters of $H_3$ give rise to bielliptic surfaces in
$|-K+T|$ where $T$ is some torsion divisor.
\skipaline

At this point we can also understand the different $G_6$-embeddings
of a scroll $X_F$ into ${\Pn 5}$. For this we start with a $G_6$-invariant
scroll $X_F$. Its desingularisation is ${\bf P}_F({\cal E})$. First note
that the Heisenberg group $H_2$ acts on the planes of $X_F$ and that this
induces actions of $H_2$ on the base curve $F$ of ${\bf P}_F({\cal E})$
and on the $3$ sections which are mapped to the singular lines of $X_F$.
To define a non-degenerate map from
${\bf P}_F({\cal E})$ to ${\Pn 5}$ is
the same as defining an isomorphism from
$H^0({\cal O}_{{\bf P}_F({\cal E})}(1))=H^0({\cal E})$ to
$V$. The
decomposition ${\cal E}={\cal L}_0\oplus h^*{\cal L}_0\oplus (h^2)^*{\cal L}_0$
defines a decomposition $H^0({\cal E})=U_1 \oplus U_2 \oplus U_3$ where
each of the spaces $U_i$ is the space of sections of a degree $2$ line bundle
on $F$ and hence has dimension $2$. The level $2$ structure on $F$ which
comes from the action of $G_6$ on $X_F$ gives us an identification
(unique up to a scalar) of
each of the $U_i$ with the $H_2$-module $U$. Now we pick one of the $4$
subgroups
$K_j$. Recall that as an $H_2$-module $V=3U$, whereas as an $H(K_j)$-module
$V=U \oplus U'\oplus {\bar {U}}'$ (cf. Lemma 7.3). In order to map the $3$
decomposing sections of ${\bf P}_F({\cal E})$ to the singular lines
associated to the group $K_j$ we must map each of the spaces $U_i$ to one
of the spaces $U$,$U'$ and ${\bar {U}}'$. The group $G_6$ acts transitively
on the lines associated to $K_j$ and hence we can (up to an element in
$G_6$) assume that $U_1$,$U_2$,$U_3$ map to $U$,$U'$, ${\bar {U}}'$. This
defines the isomorphism from $H^0({\cal E})$ to $V$ up to an element
in $({\bf C}^*)^3$. On the other hand the Heisenberg group $H_3$
acts irreducibly on the $3$-dimensional space given by the
decomposition $V=3U$. Hence, by Schur's lemma, the isomorphism from
$H^0({\cal E})$  to $V$ is uniquely defined (up to a scalar). This
shows that given a $G_6$-invariant embedding of $X_F$ in ${\Pn 5}$
we can find four such embeddings, one for each of the subgroups
$K_j$.
\skipaline
We can now consider the incidence correspondence
$${\cal I}(A,X)=\{(A,X)| A\in {\cal A}, X\in {\cal V},A \subset X\}\subset
{\cal A}\times{\cal V}.$$

\Proposition 7.10. This is a $2:3$ correspondence.

\Proof  Clearly $A\in {\cal A}$ is the intersection of two scrolls.  On
the other hand,
the number of such surfaces in a scroll is the answer to the question:
how many abelian
surfaces are there in the pencil $H\subset |-K_ {{\bf P}_F({\cal
E})}|$ coming from the Hesse pencil
which are isomorphic to the product $E\times F$? For this we want to find
an embedding
$E\subset \Pn 2\times E$ which after projection to $\Pn 2$ maps $E$ to an
element in the
Hesse pencil which is $\Z /3$-equivariant (here $\Z /3$ is the group
which acts on $\Pn 2 \times E$ with quotient ${\bf P}_F({\cal
E})$).  To embed $E$ as an element in the Hesse pencil is
the same as choosing a level 3 structure on $E$.  Say $\Z /3$ acts on $E$
by translation
with an element $\sigma^{\prime}$ of order $3$.  So we have to ask in how may
ways we can extend
$\sigma^{\prime}$ to a level 3 structure.  If $\tau^{\prime}$ is another
$3$-torsion point with
$(\sigma^{\prime},\tau^{\prime})=1$  (here $(\ , )$ is the Weil pairing), we
have the
possibilities $(\sigma^{\prime},\tau^{\prime})$,
 $(\sigma^{\prime},\tau^{\prime}\sigma^{\prime})$,
$(\sigma^{\prime},\tau^{\prime}{\sigma^{\prime}}^2)$ and no others.  This
gives us the
three possibilities.\hfill\qed

\Remark 7.11. Notice that the three choices of $\Z /3$-subgroups of $\Z
/3\times \Z /3$ correspond precisely to the three subgroups $K_i$ distinct
from the
subgroup $K_j$ which stabilizes the three singular lines of the scroll
$X_F$.

\Corollary 7.12.  Given two distinct subgroups $K_i$ and $K_j$, the
incidence
correspondence ${\cal I}(A,X)$ defines a $1:1$ correspondence between
elliptic scrolls $X$
whose singular lines are invariant under these two subgroups.

\Proof  A scroll singular along one of the triples of lines contains
exactly one
abelian surface which forms the intersection with a scroll singular along
the other
triple of lines.\hfill\qed

\Corollary 7.13.  The abelian surfaces with two plane cubic curve fibrations
are contained in
precisely a pencil of $G_6$-invariant pencils of
cubic hypersurfaces.

\Proof  The space of $G_6$-invariant pencils of cubics is a $\Pn 3$,
and the
$G_6$-invariant scrolls singular along a triple of lines are defined by
points on four lines
in $\Pn 3$ corresponding to the four subgroups $K_i$. An abelian
surface in the
intersection of two scrolls is contained in the pencils of cubics
corresponding to a
line joining two of these lines.  If there were more than a pencil of
invariant pencils of
cubics through the surface, then it would be contained in four scrolls.
This is
impossible (see Proposition 7.10). \qed

\skipaline
Consider the Grassmannian of lines in the space $\Pn 3$ of
$G_6$-invariant pencils of
cubics.   The four lines of pencils defining $G_6$-invariant scrolls are
pairwise disjoint.
The lines corresponding to abelian surfaces $A\in {\cal A}_{ij}$ define a
one to one
correspondence between two skew lines, so they form a conic section in the
Grassmannian.
Summing up we have 6 disjoint conic sections in the Grassmannian
parametrizing ${\cal A}$.

\skipaline

We can now sum up our discussion as follows:

\Proposition 7.14. The variety ${\cal A}$ consists of $6$
irreducible components ${\cal A}_{ij}$
indexed by the pairs $\{i,j\}\subset \{1,2,3,4\}$. The
incidence variety ${\cal I}(A,X)$ consists of 12 components ${\cal I}^k_{ij}$
indexed by pairs $\{i,j\}$ and an element $k\in\{i,j\}$. For every component
${\cal I}^k_{ij}$ there is a diagram
\skipaline
$$
\matrix{
X^0(6){\varphi\atop\longrightarrow} {\cal I}^k_{ij}
{q\atop\longrightarrow} {\cal V}^k\cr
{\scriptstyle p}\downarrow\cr
{\cal A}_{ij}}
$$
where $\varphi$ is $3:1$ and $p$ and $q$ are $1:1$.\par

\Proof The map $\varphi$ is the map which associates to each pair
$(E,(s_6,t_6))$ the abelian surface $A=E\times F$, the level $6$
structure $(w_1, w_2)$
and the scroll $X$ which is the $\PP^2$-scroll attached to $A$ which is
$K_k$-invariant. This map factors through ${\cal A}_{ij}$
(giving the map $\psi_{ij}$ of Lemma 7.9) and in particular the
projection $p$ has an inverse.\par
The map $q$
is $1:1$ by Corollary 7.12.\hfill\qed
\skipaline

\Remark 7.15. The
components ${\cal A}_{ij}$, ${\cal I}_{ij}^k$ and ${\cal V}^k$ are all
rational. We have already observed this for the varieties ${\cal V}^k$ which
are isomorphic to an open set of a pencil of cubics. Since the maps
$p$ and $q$ are birational, this is also true for the other varieties.
\skipaline

For our conclusion on $3$-folds $Y=X_E\cup X_F$ we study a $G_6$-invariant
rational cubic scroll.  It plays a crucial role when we later bilink $Y$ to
a $3$-fold of degree $7$ (cf. Propositions 4.7 and 5.2).\par
    The subgroups $H(K_i)$ have a
nontrivial intersection, namely the subgroup $\langle\tau^3,\sigma^3\rangle$.
This subgroup therefore fixes all 4 triples of lines.  We shall find small
$G_6$-orbits of
planes intersecting all these lines.
   Consider the subgroup $$G_3=\langle\sigma^2, \tau^2,\iota\rangle\subset
G_6.$$ We
look for
$G_3$-invariant planes.

Now the action of $\tau^2$ is defined by $$\tau^2={\rm diag}(1,\eta, \eta^2,
1,\eta, \eta^2),$$ where $\eta=e^{2\pi
i/3}$  while $\sigma^2$ sends $x_i\mapsto x_{i+2}$ and
$\iota$ sends $x_i\mapsto x_{-i}$.
  This
implies easily that the 4 planes in the $G_6$-orbit of any such plane must
be of the form

$$\matrix {
\alpha x_0+\beta x_3&=\alpha x_2+\beta x_5&=\alpha x_4+\beta x_1&=0&&P_0\cr
\alpha x_0-\beta x_3&=\alpha x_2-\beta x_5&=\alpha x_4-\beta x_1&=0&&P_1\cr
\beta x_0+\alpha x_3&=\beta x_2+\alpha x_5&=\beta x_4+\alpha x_1&=0&&Q_0\cr
\beta x_0-\alpha x_3&=\beta x_2-\alpha x_5&=\beta x_4-\alpha
x_1&=0&&Q_1.\cr }$$

In particular we notice that there is a 1-parameter family of such planes.
The union of these planes forms a rational cubic scroll $R$:  In fact the
union of the
planes is defined by
$${\rm rank}\pmatrix {x_0&x_2&x_4\cr x_3& x_5& x_1\cr}\leq 1,$$
so the scroll $R$ is the Segre embedding of $\Pn 1\times \Pn 2$ in $\Pn 5$.

\Remark 7.16.  The orbits of planes in $R$ do not all have length 4, in fact
$$P_0\cap P_1=\emptyset,\quad Q_0\cap Q_1=\emptyset$$
$${\rm and}\quad P_0\cap Q_0=P_0\cap Q_1=P_1\cap Q_0=P_1\cap
Q_1=\emptyset\Leftrightarrow \alpha^2\ne \pm\beta^2.$$
Hence the orbit is 4 disjoint planes unless $\alpha^2=\pm\beta^2$.
In case $\alpha^2=-\beta^2$ let $\alpha =1$.
Then $\beta=\pm i$ and
$P_0=Q_1, P_1=Q_0$.  While in case $\alpha^2=\beta^2$, let  $\alpha =1$.
Then $\beta=\pm 1$ and
$P_0=Q_0, P_1=Q_1$.
Thus the 4 planes are disjoint unless 2 planes come together, which happens
for
$$(\alpha,\beta)=(1,0),(0,1),(1,1),(1,-1), (1,i),(1,-i), (1,1), (1,-1),$$
altogether 4
pairs of planes.

\skipaline
\Remark 7.17. The cubic scroll $R$ contains all four sets of $G_6$-invariant
triples of lines.  These lines
are all transverse to the planes of $R$.  Thus the six singular lines of
two scrolls $X_E$
and $X_F$ which intersect along an abelian surface $E\times F$ are all
contained in a rational
cubic scroll.

\skipaline

In Lemma 8.1 we shall show that the union $Y=X_E\cup X_F$ is contained in 6
quintics. These are all singular along two
triples of lines in the scroll $R$.  We analyze these quintics more closely.
  With
bihomogeneous coordinates $s,t$ and $y_0, y_1, y_2$ on $R\cong \Pn 1\times
\Pn 2$, the
restriction of quintics singular along the two triples of lines generated
by $\tau$ and
$\sigma$ have the form  $$\eqalign
{y_0^2y_1^3-y_0^3y_1y_2-y_0y_1^2y_2^2+y_0^2y_2^3\cr
y_0^3y_1^2-y_0y_1^3y_2-y_0^2y_1y_2^2+y_1^2y_2^3\cr
y_0^2y_1^2y_2-y_0^3y_2^2-y_1^3y_2^2+y_0y_1y_2^3,\cr}$$
multiplied by any quintic in $s,t$.
The 6 quintics in the ideal of two scrolls are determined by a
pencil of quintics in $s,t$.
This pencil is $G_3$-invariant, so it is an element in the net
$$\langle s^5,t^5\rangle\oplus \langle s^4t,st^4\rangle\oplus \langle
s^3t^2, s^2t^3\rangle$$ of
pencils. The basepoints of this
pencil define precisely the planes of $R$ common to all quintics through $W$.
 If the pencil has no basepoints, then the two scrolls would intersect
$R$ in only the
singular lines.  We shall see in Proposition 7.18 that this
is not the case.  On the
other hand every point of intersection of the two scrolls with $R$
outside the 6 singular lines lies in a
plane which must be defined by a basepoint of the pencil, so there are at
most 4 planes with
such an intersection. These planes clearly form orbits under $G_6$ so there
are 2 or 4 planes
as explained above.

\skipaline

We carry this analysis a bit further in order to show that the union of two
scrolls is bilinked to a Del Pezzo $3$-fold $W_7$.
First we consider  the intersection of the cubic
scroll $R$ and an elliptic scroll $X$.  We use the fact that $X$ is
$G_6$-invariant and that each plane in $R$ is invariant  under the subgroup
$G_3\subset G_6$ generated by $\langle\sigma^2, \tau^2\rangle$.  Thus any
point in a
plane has an orbit by
this subgroup of order divisible by $3$, and any invariant curve
 has
degree divisible by $3$.

\Proposition 7.18. The intersection $X\cap R$ is a curve of degree $18$.
It has
the following decomposition into irreducible components:
 $$C=2L_1+2L_2+2L_3+l_1+\ldots +l_{12}$$
where $L_i$ are the singular lines of $X$ and as such meet every
plane in $R$.  The lines
$l_i$ form $4$ triangles in $4$ planes of the $\Pn 2$-bundle $R$.

\Proof  First we prove that the intersection is a curve.  If not,
 it
contains an irreducible
surface, call it $T$.  If $T$ contains a plane, this is common to
$R$ and $X$.  But
no plane in $X$ is stabilized by $G_3$, in fact $X$ intersects only one of
the 4 triples of
lines, so this is impossible. Thus $T$ intersects each plane in $R$ in some
curve.  This
curve has degree $3$ since $X$ is contained in two cubics and the curve is
invariant under
$G_3$.

As noted above, the three singular lines of the
scroll $X$ all lie in $R$, and are transverse to the planes in $R$.
Since every cubic through $X$ is singular along the three lines, the
intersection of a
cubic with $R$ is a triangle in each plane.  In fact $T$ in $R=\Pn 2\times
\Pn 1$ must
equal $T=T_0\times \Pn 1$, where $T_0$ is this triangle, i.e. $T$ is the
union of three
quadric surfaces.  But $X$ does not contain any quadric surface, so the
intersection
$X\cap R$ cannot contain a surface. It is therefore a curve, call it
$\Gamma$.  Since
the intersection is proper, this curve has, by Bezout, degree $18$.
\skipaline
We have seen
already that the curve  $\Gamma$ contains the three singular lines. In fact on
$\tilde X$, the normalization of $X$, the preimage $\tilde\Gamma$ of $\Gamma$
contains the curves $F_i$ which lie $2:1$ over the singular lines $L_i$.
Thus on $\tilde X$ we have
$$\tilde\Gamma =F_1+F_2+F_3+ah^2+bhf,$$
in notation as in the proof of Proposition 4.10, with $6a+b=12$.  Since any
plane which
intersects $R$ properly intersects in a scheme
of length 3, the general plane of $X$ must intersect $R$ in three points,
the points of
intersection between the plane and the lines $L_i$.  Therefore $a=0$, and
$\tilde\Gamma =F_1+F_2+F_3+l_1+\ldots +l_{12}$ where the $l_i$ are lines in
the planes
of $\tilde X$.

Finally if a plane in $R$ intersects $X$ along a curve, this curve has
degree divisible
by three, in fact equal three, since $X$ lies in two cubics.  Therefore the
twelve lines $l_i$
form triangles in four planes of $R$. \qed

\Section 8. Conclusion

We shall conclude by showing that the union of two scrolls $X_E$ and
$X_F$ which intersect along an abelian surface $E\times F$ is
bilinked to a Del Pezzo $W_7$.

Let $Y=X_E\cup X_F$.

\Lemma 8.1. $\cohi 0Y5=6$ and $Y$ lies on irreducible quintic
hypersurfaces.\par

\Proof
 Consider the
exact sequence
$$0\to {\id Y5}\to{\id {X_E}5}\oplus{\id {X_F}5}\to {\id {E\times F}5}\to 0.$$

First, the intersection  $E\times F\cap\Pn 3$ of $E\times F$ with a general
$\Pn 3$ is at least $5$-regular in the sense of of Castelnuovo-Mumford, so
the same
is true for $E\times F$, i.e. $\cohi 1{E\times F}k=0$, when $k\geq 4$.
In particular,  $\cohi 0{E\times F}5=102$. Similarly,
the intersection $Y\cap \Pn 2$ of $Y$ with a general plane is $5$-regular
in the
sense of Castelnuovo-Mumford so $\cohi 1Yk=0$ for $k\geq 4$. Furthermore $\cohi
0{X_E}5=\cohi 0{X_F}5=54$ by Corollary 6.4. Therefore $\cohi 0Y5=6$. If every
quintic in the ideal of $Y$ is reducible, then they all have a fixed quartic
hypersurface as a component.  But $Y$ is $G_6$-invariant, and there are no
$G_6$-invariant quartic hypersurfaces, so the lemma follows.\qed

Now, according to Proposition 7.18 the two elliptic scrolls $X_E$ and $X_F$
intersect the rational scroll $R$ in the six singular lines and in
four triangles each.  Since $X_E\cap X_F$ is a surface which does not meet
the six singular
lines,
the two sets of four triangles must lie pairwise in the same four planes of
$R$.  Thus the six
quintics through the two scrolls $X_E\cup X_F$ contain these four planes.
In particular,
if $Z$ is linked $(5,5)$ to the union $X_E\cup X_F$, then $Z$ has
degree $13$ and contains four planes of $R$.
Furthermore since each quintic intersects the planes  in $R$ in curves
singular in the six
points of intersection with the singular lines, $Z$ will intersect each
plane in one point
in addition to the six singular points.  The general quintic in the pencil
intersects a plane in an irreducible curve, so, by Bezout, there are no
quartic curves singular in the six singular points.  Therefore any quartic
hypersurface containing $Z$ must also contain $R$.

Now, the partial normalization
$Y^{\prime}$ of $Y=X_E\cup X_F$ along the $6$ singular lines is Calabi-Yau (cf.
Proposition 3.5). In particular the dualizing sheaf $\omega_{Y^{\prime}}$
is trivial and has one
global section. It follows from the next lemma that the dualizing sheaf
$\omega_Y$ also has a
section.

Let $f:Y^{\prime}\to Y$ be a finite morphism of projective schemes and
denote by $\omega_Y$ the
dualizing sheaf of $Y$.  Then $f^!\omega_Y$ is a dualizing sheaf of
$Y^{\prime}$ (See [Ha, Ex
III, 7.2] and for the definition of $f^!\omega_Y$ see [Ha, Ex III, 6.10]).
Hence we can put
$\omega_{Y^{\prime}}=f^!\omega_Y$.

\Lemma 8.2.  If ${\rm H}^0({Y^{\prime}},\omega_{Y^{\prime}})\not=0$, then
also ${\rm
H}^0(Y,\omega_Y)\not=0$.\par

\Proof  Using [Ha, Ex III,6.10(b)] we have that
$$\eqalign {{\rm H}^0({Y^{\prime}},\omega_{Y^{\prime}})&={\rm Hom}_{{\cal
O}_{Y^{\prime}}}
({\cal O}_{Y^{\prime}},\omega_{Y^{\prime}})= {\rm Hom}_{{\cal O}_{Y^{\prime}}}
({\cal O}_{Y^{\prime}},f^!\omega_Y)\cr & = {\rm H}^0({{\cal H}om}_{{\cal
O}_{Y^{\prime}}} ({\cal O}_{Y^{\prime}},f^!\omega_Y))\cr
& = {\rm H}^0(f_*{{\cal H}om}_{{\cal O}_{Y^{\prime}}}
({\cal O}_{Y^{\prime}},f^!\omega_Y))\cr
& = {\rm H}^0({{\cal H}om}_{{\cal O}_Y}
(f_*{\cal O}_{Y^{\prime}},\omega_Y)) \quad ([Ha., Ex. III. 6.10(b)])\cr
& = {\rm Hom}_{{\cal O}_Y}
(f_*{\cal O}_{Y^{\prime}},\omega_Y).}$$
Hence a section $s$ of $\omega_{Y^{\prime}}$ gives rise to a morphism
$\varphi_s:f_*{\cal O}_{Y^{\prime}}\to\omega_Y$.
Combining this with the natural morphism ${\cal O}_Y\to f_*{\cal
O}_{Y^{\prime}}$
 we obtain a morphism ${\cal
O}_Y\to \omega_Y$ and hence a section of $\omega_Y$.\qed
\skipaline

  In the cohomology of the liaison exact sequence (cf. [PS]) $$0\to
\omega_Y\to {\cal O}_{Y\cup Z}(4)\to {\cal O}_{Z}(4)\to 0,$$ a section of
$\omega_Y$ corresponds to a section of $\cohi 0Z4$, i.e. to a  quartic
hypersurface
containing $Z$. This quartic must contain $R$.

 Let  $W$ be linked to $Z$ in this
quartic and a general quintic through $Z$.  Then $W$ must intersect $R$ in a
surface, in fact in a surface linked to 4 planes in the intersection of $R$
with a
quintic hypersurface.  This is clearly a surface of degree $11$.
Furthermore the
arithmetic genus of $W$ is 1 by linkage, so
$W\cup R$ has arithmetic genus 11.  From the liaison exact sequences

$$0\to \omega_Z(1)\to {\cal O}_{Y\cup Z}(5)\to
{\cal O}_{Y}(5)\to 0,$$
$$0\to \omega_{Z}(1)\to {\cal O}_{W\cup Z}(4)\to
{\cal O}_{W}(4)\to 0$$
we get that $\cohi 0W4=\cohi 0Y5-1=5$.

Clearly all the quartics through $W$ have to contain also
$R$.  Thus $W\cup R$ is contained in 5 quartics.
\par
 For later we need that the 5 quarties are
minors of a $4\times 5$ matrix with linear entries.  If we can show this in
a special case, we may
conclude by semicontinuity that
so is $W\cup R$ in general (cf. [Ell]).

We do this by considering, in the notation of Remark 7.8,  the
points $(1,-3,0,0)$ and $(0,0,0,1)$ in the $\Pn 3$ of
$G_6$-invariant pencils of cubics.  They are points on two distinct
lines corresponding to two distinct subgroups $H(K_i)$.  The
corresponding pencils of cubics each define an elliptic
scroll (in fact a reducible scroll) residual to three $\Pn 3$'s. The union
of these scrolls
lies on 6 quintics and is bilinked $(5,5)$ and $(5,4)$ to a
$3$-fold of degree $7$ which lies in 5
quartic hypersurfaces.  These 5 quartics define a determinantal $3$-fold of
degree $10$. These claims are easily checked
with i.e. [MAC].
\par
Recall from Proposition 4.10 that $W$ is a non-normal Del-Pezzo $3$-fold $W_7$
as soon as we have checked that no three of the singular lines lie in a
$\Pn 3$,
and that not all six lie in a rational normal quartic scroll, and finally
that the common
singular locus of the quartic hypersurfaces through $W$ is precisely the
$6$ singular
lines.  While the former two requirements follows easily from our analysis
of the singular
lines in the previous section, the latter requirement is easily checked in
the above
example.  Therefore we may, by Proposition 5.2, bilink
$W$ in complete intersections $(5,4)$ and $(5,5)$ to a non-normal
Calabi-Yau $3$-fold $Y_t$, non-normal only along 6 lines.  Clearly, we may
perform the
bilinkage in a family, so the reducible $3$-fold $Y$ is a degeneration
of $Y_t$. We have shown

\Theorem 8.3.  The reducible $3$-fold $Y=X_E\cup X_F$
is a degeneration of irreducible non-normal Calabi-Yau $3$-folds of degree $12$
in $\Pn 5$. The general such $3$-fold is singular precisely along $6$ disjoint
lines.\par

\Remark 8.4.  Computing the normalizer $N_3$ of $G_3\subset GL(6,{\bf C})$
and its
representations, one may show that there is an $N_3$-invariant linear
complex in the
Grassmannian $G$ of lines in the space of $G_6$-invariant pencils of
cubics. This defines a
one to one correspondence between any two of the four lines defined by the
subgroups $K_i$.
Therefore it is natural to guess that this complex defines the
correspondences ${\cal A}_{ij}$, in
fact that it parametrizes the set of all $G_6$-invariant abelian surfaces.
This is proved
by Gross and Popescu [GP2].
\skipaline

\References MAC

\ref Ch Chang, M.-C., Classification of Buchsbaum subvarieties of
codimension $2$ in
projective space.  J. Reine Angew. Math. 401 (1989), 101--112.

\ref DP Decker, W.,  Popescu, S., On surfaces in $\Pn 4$ and $3$-folds in
$\Pn 5$. In:
Vector bundles in algebraic geometry, eds. Hitchin, N.J., Newstead, P.E.,
Oxbury, W.M.
Cambridge University Press L.N.S. 208 (1995), 69--100.

\ref Ell Ellingsrud, G.,  Sur le schema de Hilbert des vari\'et\'es de
codimension 2  dans $\Pn e$ \`a cone de Cohen Macaulay. Ann. Scient. E.N.S. 4.
Serie, 8 (1975), 423--432.

\ref ES Ellingsrud, G., Str\o mme, S. A.,  The number of twisted cubic
curves on the general quintic threefold. Math. Scand. 76 (1995), 5--34.

\ref Fu Fulton, W.,  Intersection Theory. Springer Verlag, New
York-Berlin-Heidelberg (1982).

\ref GP1  Gross, M.,  Popescu, S.,  Equations of
$(1,d)$-polarized abelian surfaces.  Math. Ann. 310 (1998), 333--377.

\ref GP2  Gross, M.,  Popescu, S.,  Calabi-Yau $3$-folds and moduli of abelian
surfaces, I and II.  To appear.

\ref Ha  Hartshorne, R., Algebraic Geometry. GTM {\bf 52} Springer-Verlag,
New York,
1977.

\ref Kre Kresch, A., FARSTA, a computer program for quantum cohomology.
Appendix to: Quantum Cohomology at the Mittag-Leffler Institute,
 edited by P. Aluffi, Scuola Normale Superiore, Pisa (1997).

 \ref MAC Bayer, D., Stillman, M.,  MACAULAY: A system for
computation in algebraic geometry and commutative algebra, Source and
object code available for Unix and Macintosh computers.  Contact the
authors, or download from zariski.harvard.edu via anonymous ftp.

\ref PS Peskine, C., Szpiro, L.,  Liaison des vari\'et\'es alg\'ebriques I.
Invent. Math.  26 (1974), 271--302.

\ref Re Reider, I., Vector bundles of rank $2$ and linear systems on algebraic
surfaces. Annals of Math. 127 (1988), 309--316.

\ref Rey Reye, T., Ueber lineare Systeme und Gewebe von Fl\"achen
zweiten Grades. J. Reine Angew. Math. 82 (1877), 54--83.

\skipaline
Authors' email addresses:
  hulek@math.uni-hannover.de, ranestad@math.uio.no

  \bye